\documentclass[12pt]{amsart}
\usepackage{amsmath, amsthm, amssymb}
\usepackage[margin=1in]{geometry}
\usepackage{epsfig}
\usepackage{rawfonts}
\usepackage{enumerate}
\usepackage{graphics}
\usepackage{multirow}
\usepackage{xspace}
\usepackage{graphicx}
\usepackage{pgf,tikz,pgfplots}
\usepackage{mathrsfs}
\usetikzlibrary{arrows}
\usepackage{amsmath}
\usepackage{amsfonts}
\usepackage{amssymb}
\usepackage{amsthm}
\usepackage{graphicx}
\usepackage{booktabs}
\usepackage{caption}
\usepackage{listings}
\usepackage{setspace}
\usepackage[mathscr]{eucal}
\usepackage{pgfplots}
\usepackage{hyperref}
\usepackage{wrapfig}
\usepackage{floatflt,epsfig}
\usepackage{ dsfont }
\usepackage{amscd}
\usepackage{tikz-cd}
\usepackage{fancyhdr}
\usepackage[all]{xy}
\usepackage{latexsym}
\usepackage{amscd}
\usepackage{pifont}
\usepackage{eufrak}
\usepackage{subfig}
\usepackage{easyReview}

\newcommand{\numberset}{\mathbb}
\newcommand{\N}{\numberset{N}}

\theoremstyle{plain}

\theoremstyle{theorem}

\newtheorem{definition}{Definition}[section]
\newtheorem{proposition}[definition]{Proposition}
\newtheorem{theorem}[definition]{Theorem}
\newtheorem{lemma}[definition]{Lemma}
\newtheorem{conjecture}[definition]{Conjecture}
\newtheorem{corollary}[definition]{Corollary}
\newtheorem{example}[definition]{Example}
\newtheorem{remark}[definition]{Remark}

\theoremstyle{remark}

\makeatletter
\@namedef{subjclassname@2020}{\textup{2020} Mathematics Subject Classification}
\makeatother

\begin{document}
		
		\title[Structure of gaps and generalized Wilf's conjecture]{Structure of gaps and generalized Wilf's conjecture}
		
			\author{OM PRAKASH BHARDWAJ}
		\address{Chennai Mathematical Institute, Tamil Nadu - 603103, India}
		\email{omprakash@cmi.ac.in, opbhardwaj95@gmail.com}

		\author{CARMELO CISTO}
		\address{Universit\`{a} di Messina, Dipartimento di Scienze Matematiche e Informatiche, Scienze Fisiche e Scienze della Terra\\
			Viale Ferdinando Stagno D'Alcontres 31\\
			98166 Messina, Italy}
		\email{carmelo.cisto@unime.it, cistomelo@gmail.com}
		
		\author{Priyavrat Deshpande}
		\address{Chennai Mathematical Institute, Tamil Nadu - 603103, India}
		\email{pdeshpande@cmi.ac.in}

		
		\begin{abstract}
 In this article, we study the gaps, the maximal gaps, and the pseudo-Frobenius elements of generalized numerical semigroups of minimum embedding dimension. For these semigroups, we give an explicit description of the set of gaps, explicit formulas for the genus and the number of maximal gaps, and an explicit description of the set of pseudo-Frobenius elements in terms of the blocks of consecutive gaps of an associated numerical semigroup. In particular, we obtain an explicit formula for the type. As a consequence, we show that the type is unbounded even in minimum embedding dimension. Further, we prove that every generalized numerical semigroup with minimum embedding dimension has at least as many non-gaps below its gaps as gaps, and therefore satisfies the generalized Wilf's conjecture. We also prove the conjecture for generalized numerical semigroups in $\mathbb{N}^2$ graded by two-generated numerical semigroups, and asymptotically in higher dimension. For an arbitrary generalized numerical semigroup we bound the number of non-gaps below the gaps in terms of the genus and the type. Finally, we introduce nearly-symmetric generalized numerical semigroups and classify them in minimum embedding dimension. Although the type is unbounded in minimum embedding dimension generalized numerical semigroups, nearly-symmetric ones have type at most two.
  	\end{abstract}

		\maketitle
		
\section{Introduction}
Affine semigroups and semigroup rings have a rich literature and connections to combinatorics, commutative algebra, and algebraic geometry. There exists a one-to-one correspondence between normal affine semigroups rings and and affine toric varieties. In the one-dimensional case, affine semigroup rings are isomorphic to numerical semigroup rings. Numerical semigroups play a crucial role in the study of one-dimensional local domains. Moreover, numerical semigroups naturally arise in the study of algebraic curves. Let $P$ be a rational point on an algebraic curve $C$. 
The set of integers $k$ for which there exists a function on $C$ having pole divisor exactly $kP$ forms a numerical semigroup, known as the Weierstrass semigroup of $P$.  
Furthermore, when $P_1, \ldots, P_d$ are distinct points on $C$, 
one defines the Weierstrass gapset as

{\Small
	\[
	G(P_1, \ldots, P_d) := 
	\big\{ (y_1, \ldots, y_d) \in \mathbb{N}^d 
	\ \big|\ 
	\nexists \text{ a meromorphic function on } C 
	\text{ with pole divisor } y_1 P_1 + \dots + y_d P_d \big\},\]}
    
    then $G(P_1, \ldots, P_d)$ is finite and $\mathbb{N}^d \setminus 	G(P_1, \ldots, P_d)$ is an affine semigroup. These semigroups and their gapsets have been studied extensively in the literature, for example, see (\cite{Ballico-Kim98,Beelen2006,Carvalho-et-al2005,Ishii2002,Matthews2001}). In general, an affine semigroup in $\mathbb{N}^d$ which has finite complement in $\mathbb{N}^d$ is called a \textit{generalized numerical semigroup}.

Generalized numerical semigroups were formally introduced by Failla, Peterson and Utano~\cite{Failla-et-al2016}, and have since been studied from many points of view, for example, see (\cite{Bernardini-et-al2024,Bhardwaj-Goel-Sengupta2023, Cisto-Tapia-Ramos2026,Failla-et-al2016, Garcia-et-al2024, Li2025,HemaSrinivasan2025,Singhal-Lin2022}).  Gaps, genus, embedding dimension, pseudo-Frobenius elements and type extend naturally to this setting, but the total order of $\mathbb{N}$ does not extend to $\mathbb{N}^d$. Therefore, in general, a generalized numerical semigroup does not have a largest gap, and hence no canonical Frobenius element. To overcome this situation, one either fixes a monomial order and take the maximum element of the gap set~(\cite{garcia2018extension, garcia-pseudofrobenius}) or studies the whole
poset of gap set and its maximal elements, called Frobenius allowable gaps~(\cite{Singhal-Lin2022}).  In this article we take the second point of view. We describe the gaps, the maximal
gaps and the pseudo-Frobenius elements of a natural extremal family explicitly (generalized numerical semigroups with minimum embedding dimension),
and use this description to address generalized Wilf’s conjecture, type, the relation between gaps and non-gaps, and nearly-symmetric property.

It is proved in \cite{garcia2018extension} that every generalized numerical semigroup $S \subsetneq \mathbb{N}^d$ has at least $2d$ minimal
generators. Therefore the generalized numerical
semigroups with embedding dimension $2d$
are extremal objects in this theory, playing
the role that two-generated numerical semigroups
play for $d=1$. We call these generalized
numerical semigroups the minimum embedding
dimension generalized numerical semigroups.
Define the following generalized numerical
semigroup:
$$\mathcal{S}=\mathcal{S}(T,\{d_j\}_{j \neq i})
=\left\langle
\{a_1\mathbf{e}_i,a_2\mathbf{e}_i,\ldots,a_n\mathbf{e}_i\}
\cup
\left\{
\mathbf{e}_j, \mathbf{e}_i +d_j\mathbf{e}_j
\mid j\in\{1,\ldots,d\}\setminus \{i\}
\right\}
\right\rangle,$$
where $T=\langle a_1,a_2,\ldots,a_n\rangle$ is a numerical semigroup minimally generated by $\{a_1,a_2,\ldots,a_n\}$, $d_j$ are positive integers, and $\mathbf{e}_1,\ldots,\mathbf{e}_d$ are the standard basis vectors of the vector space $\mathbb{R}^d$. When
$T=\langle a_1,a_2\rangle$
then $\mathcal{S}$ is minimally generated by $2d$ elements
and all generalized numerical semigroups of
minimal embedding dimension are of this
form up to a permutation of coordinates. Though the class of generalized numerical semigroups
with minimum embedding dimension is already
a very important extremal class, in this
article, we explore a more general class
of generalized numerical semigroups $\mathcal{S}$ as defined above.

Although every generalized numerical semigroup of minimum embedding dimension is of this form with \(T\) two-generated, we allow \(T\) to be an arbitrary numerical semigroup. 
This broader family allows us to describe the gaps and maximal gaps in terms of the consecutive blocks of gaps of \(T\). 
In particular, the genus and the number of Frobenius allowable gaps can be expressed explicitly in terms of these blocks.

For a gap $h$ of the numerical semigroup $T$, let $\alpha(h)$ denote the number of consecutive gaps of $T$ ending at $h$. We first describe the gaps of $\mathcal{S}=\mathcal{S}(T,\{d_j\}_{j\neq i})$ in terms of the gaps of $T$. More precisely, we prove that an element $\mathbf{x}\in\mathbb{N}^d$ is a gap of $\mathcal{S}$ if and only if $x_i$ is a gap of $T$ and
$$\sum_{j\neq i}\left\lfloor\frac{x_j}{d_j}\right\rfloor
\leq \alpha(x_i)-1.$$
 Using this description, we obtain an explicit formula for the genus of $\mathcal{S}$. If the gaps of $T$ are decomposed into maximal blocks of consecutive integers of lengths $\ell_1,\ldots,\ell_r$, then
$$\operatorname{g}(\mathcal{S})=\Big(\prod_{j\neq i}d_j\Big)
\sum_{\nu=1}^{r}\binom{\ell_\nu+d-1}{d}.$$
In particular, the genus of $\mathcal{S}$ depends on $T$. Furthermore, we give a complete description of the maximal gaps of $\mathcal{S}$ and obtain a formula for the number of maximal gaps (Theorem~\ref{thm:Gaps}, Theorems~\ref{thm:genus} and~\ref{thm:FA}).

Next, we determine the type and the pseudo-Frobenius elements of the generalized numerical semigroup $\mathcal{S}$. The pseudo-Frobenius elements of an affine semigroup $S$ are closely related to the homological properties of its semigroup ring. More precisely, they determine the degrees of the last syzygies of the associated semigroup ring $\mathbb{K}[S]$, and the number of pseudo-Frobenius elements, called the type $\operatorname{t}(S)$, is the corresponding Betti number~(see~\cite{garcia-pseudofrobenius}). We call a maximal block $I$ of consecutive gaps of $T$ admissible if $I+a_l\not\subseteq \mathcal{H}(T)$ for every minimal generator $a_l$ of $T$. We prove that the pseudo-Frobenius elements of $\mathcal{S}$ are precisely the elements of the form
$$h\mathbf{e}_i+\sum_{j\neq i}(v_jd_j-1)\mathbf{e}_j,$$ 
where $h$ is the largest element of an admissible block $I$ of gaps of $T$, and $v_j$ are positive integers satisfying $\sum_{j\neq i}v_j=|I|+d-2.$ Consequently, we obtain the following formula for the type of $\mathcal{S}$:
$$\operatorname{t}(\mathcal{S})=\sum_{I\text{ admissible}}
\binom{|I|+d-3}{d-2}.$$
In particular, when $T=\langle a,a+1\rangle$, every maximal block of gaps is admissible and
$$ \operatorname{t}(\mathcal{S})=\binom{a+d-3}{d-1}.$$
Thus, the type of a generalized numerical semigroup with minimum embedding dimension is not bounded in terms of its embedding dimension (Theorem~\ref{thm:PF} and Corollary~\ref{cor:type}).

In 1978, Wilf~\cite{Wilf1978} conjectured that for every numerical semigroup $S$, we have 
$$\operatorname{F}(S)+1\leq \operatorname{e}(S)\cdot
|\{s\in S\mid s<\operatorname{F}(S)\}|,$$
where $\operatorname{F}(S)$ and $\operatorname{e}(S)$ denote the Frobenius number and the embedding dimension of $S$, respectively. This conjecture has been studied extensively in the literature and is still open, although it has been verified for several classes of numerical semigroups; see \cite{DelgadoSurvey} for a survey. In higher dimensions, there is no Frobenius number, and the natural analog of the set $\{s\in S\mid s< \operatorname{F}(S)\}$ is the set $\operatorname{N}(S)$ of elements of $S$ lying below some gap. In~\cite{Cisto-et-al2020}, Cisto et al. proposed the generalized Wilf's conjecture. The generalized Wilf's conjecture states that for a generalized numerical semigroup $S$, we have 
$$(\operatorname{e}(S)-d)\operatorname{n}(S)\geq d \operatorname{g}(S),$$
where $\operatorname{n}(S)=|\operatorname{N}(S)|.$
For $d=1$, this is  Wilf's conjecture. For numerical semigroups, one always has $\operatorname{n}(S)+\operatorname{g}(S)=\operatorname{F}(S)+1$ and $\operatorname{n}(S)\leq \operatorname{g}(S)$. The equality holds if and only if $S$ is symmetric. In the case of generalized numerical semigroups, $\operatorname{n}(S)$ can exceed $\operatorname{g}(S)$. We prove that this is the case for the class of generalized numerical semigroups of the form $\mathcal{S}$. In particular, we construct an explicit injection from the set of gaps of $\mathcal{S}$ into $\operatorname{N}(\mathcal{S})$, and consequently prove that $\operatorname{g}(\mathcal{S})\leq \operatorname{n}(\mathcal{S}).$ As a consequence, it follows that every generalized numerical semigroup with minimum embedding dimension satisfies the generalized Wilf's conjecture (Theorem~\ref{thm_wilf-gen} and Corollary~\ref{cor:wilfminembdim}).

Next, we study the generalized Wilf's conjecture for the generalized numerical semigroups graded by numerical semigroups. These semigroups are defined as
 $$ S^d_T=\{\mathbf{x}\in\mathbb{N}^d\mid x_1+\cdots+x_d\in T\}.$$
When $T=\langle a,b\rangle$, we give explicit formulas for the genus and $\operatorname{n}(S^2_T)$, and prove that $S^2_T$ satisfies the generalized Wilf's conjecture. Furthermore, we prove that $S^d_T$ satisfies the generalized Wilf's conjecture for all sufficiently large $b$ (Theorems~\ref{thm:graded-d2} and~\ref{thm:graded-d}).

For a numerical semigroup, the number of gaps and the number of non-gaps below the Frobenius number are closely related. However, this is not true for generalized numerical semigroups. In this article, we establish a relation between the number of gaps, the number of non-gaps and the type of an arbitrary generalized numerical semigroup. More precisely, if $\tau(S)$ denotes the number of Frobenius allowable gaps of $S$, then we prove that
\[
\frac{\operatorname{g}(S)}{\operatorname{t}(S)}
\leq \operatorname{n}(S)
\leq \tau(S)\operatorname{g}(S)
\leq \operatorname{t}(S)\operatorname{g}(S).
\]
Using these bounds, we prove that for any family of generalized numerical semigroups with bounded type, the invariants $\operatorname{n}(S)$ and $\operatorname{g}(S)$ grow at the same rate (Theorem~\ref{thm:sandwich} and Corollary~\ref{cor:growth}). 

In~\cite{Herzog-et-al2019}, Herzog, Hibi and Stamate introduced the notion of nearly Gorenstein rings. Later, Moscariello and Strazzanti \cite{Moscaeiello-Strazzanti2021} characterized nearly Gorenstein numerical semigroup rings in terms of their pseudo-Frobenius elements. The study of these semigroups is important to the problem of bounding the type of a numerical semigroup in terms of its embedding dimension. This problem is still open in general, although positive results are known for almost symmetric numerical semigroups of embedding dimension less than six \cite{Moscariello16,Moscariello23} and for nearly Gorenstein numerical semigroups of embedding dimension at most four \cite{Moscaeiello-Strazzanti2021}. For $d\geq 2$, the semigroup ring associated to a generalized numerical semigroup in $\mathbb{N}^d$ is never Cohen--Macaulay, therefore, we use the pseudo-Frobenius condition to define the notion of nearly-symmetric generalized numerical semigroups. We completely characterize nearly-symmetric generalized numerical semigroups with minimum embedding dimension. In particular, we prove that a generalized numerical semigroup in $\mathbb{N}^d$ with  embedding dimension $2d$ is nearly-symmetric if and only if it is symmetric, or $d=2$ and
$$S=\langle(3,0),(4,0),(0,1),(1,1)\rangle,$$
up to permutation of coordinates. Therefore, we prove that the type of a nearly-symmetric generalized numerical semigroup with minimum embedding dimension is at most two~(Corollaries~\ref{cor:NS-classification} and~\ref{cor:NS-type}).

The article is organized as follows. In Section~\ref{sec:prelim}, we recall some preliminaries. In Section~\ref{sec:gaps}, we describe the gaps, the genus, and the maximal gaps of the semigroups of the form $\mathcal{S}:= \mathcal{S}(T,\{d_j\})$. In Section~\ref{sec:PF}, we compute their pseudo-Frobenius elements and type. Section~\ref{sec:Wilf} is devoted to the generalized Wilf's conjecture. In Section~\ref{sec:gaps-nongaps}, we study the relationship between the number of gaps and non-gaps of an arbitrary generalized numerical semigroup. Finally, in Section~\ref{sec:NS}, we introduce and explore  nearly-symmetric generalized numerical semigroups.

\section*{Acknowledgements}
The first author is supported by ANRF National Postdoctoral Fellowship. 
The first and the third authors are partially supported by a grant from the Infosys Foundation. The second author is a member of INdAM-GNSAGA. We acknowledge the use of ChatGPT Go for proofreading and help with the Sage code.

\section{Preliminaries}\label{sec:prelim}

Let $\mathbb{N}$ be the set of nonnegative integers and $d$ be a positive integer. In the set $\mathbb{N}^d$ we define $\leq_{\mathbb{N}^d}$ to be the natural partial order, that is, for $\mathbf{x},\mathbf{y}\in \mathbb{N}^d$ we have $\mathbf{x}\leq_{\mathbb{N}^d} \mathbf{y}$ if $\mathbf{y}-\mathbf{x} \in \mathbb{N}^d$.

\begin{definition}
	{\rm An additive submonoid $S$ of $\mathbb{N}^d$ is called  a \textit{generalized numerical semigroup} if the set $\mathcal{H}(S)=\mathbb{N}^d\setminus S$ is finite. In particular, when $d=1$, $S$ is called a \textit{numerical semigroup}.}
\end{definition}

  The set $\mathcal{H}(S)$ is called the set of \textit{gaps} of $S$. The cardinality of the set $\mathcal{H}(S)$ is called the \textit{genus} of $S$, and it is denoted by $\operatorname{g}(S)$. When $S$ is a numerical semigroup, the largest integer belonging to $\mathcal{H}(S)$ is called the Frobenius number of the numerical semigroup $S$, denoted by $\mathrm{F}(S)$.  
  
  Let $A$ be a subset of $\mathbb{N}^d$, we denote $$\langle A \rangle=\{\sum_{i=1}^n \lambda_i \mathbf{a}_i\mid n \in \mathbb{N}\setminus \{0\}, \lambda_i \in \mathbb{N}\setminus \{0\}, \mathbf{a}_i \in A, \text{ for all } i\in \{1,\ldots,n\}\},$$ that is the submonoid of $\mathbb{N}^d$ generated by $A$. It is known that for every generalized numerical semigroup $S$ there exists a unique finite minimal set of generators of it, that is, there exists a unique finite set $A$ such that $S=\langle A \rangle$ and for all $B\subsetneq A$ we have $S\neq \langle B \rangle$. The cardinality of the minimal set of generators of $S$ is called the \textit{embedding dimension} of $S$ and denoted by $\operatorname{e}(S)$. As shown in \cite[Theorem 11]{garcia2018extension}, the embedding dimension of a generalized numerical semigroup has the following lower bound:

\begin{proposition}\cite{garcia2018extension}
	Let $S\subsetneq \mathbb{N}^d$ be a generalized numerical semigroup. Then $\operatorname{e}(S)\geq 2d$.
\end{proposition}

\begin{definition}
 {\rm   We call a generalized numerical semigroup $S \subsetneq \mathbb{N}^d$ of\textit{ minimum embedding dimension} if $\operatorname{e}(S) = 2d$.}
\end{definition}

\begin{definition}
	{\rm Let $S\subseteq \mathbb{N}^d$ be a generalized numerical semigroup. An element $\mathbf{h} \in \mathcal{H}(S)$ is called a \textit{pseudo-Frobenius element} if for all non-zero elements $\mathbf{s} \in S$, $\mathbf{h}+\mathbf{s} \in S$. The set of all pseudo-Frobenius elements of $S$ is denoted by $\mathrm{PF}(S).$ The cardinality of the set $\mathrm{PF}(S)$ is called the \textit{type} of $S$, and denoted by $\operatorname{t}(S).$}
\end{definition}

\begin{remark}
   {\rm Let $\mathbb{K}$ be a field and $\mathbb{K}[S]$ be the semigroup ring of $S$. If $\mathbf{a}_1, \ldots, \mathbf{a}_n$ are the minimal generators of $S$ then $\mathbf{b} \in S$ is the minimal degree of the last syzygy of $\mathbb{K}[S]$ if and only if 
   $\mathbf{b} \in \{\mathbf{f} + \sum_{i = 1}^n \mathbf{a}_i \mid \mathbf{f} \in \mathrm{PF}(S)\} $~(\cite[Corollary~7]{garcia-pseudofrobenius}), and $\operatorname{t}(S)$ is the last betti number of $\mathbb{K}[S]$ which is equal to the type $\operatorname{t}(\mathbb{K}[S])$ of the semigroup ring $\mathbb{K}[S]$. From algebraic point of view, these facts motivate the notions studied in this article. }
\end{remark}

 For $\mathbf{x}, \mathbf{y} \in S$, we write $\mathbf{x} \leq_{S} \mathbf{y}$ if $\mathbf{y} - \mathbf{x} \in S$. Note that $\leq_{S}$ is a partial order on $\mathbb{N}^d$, and $\mathbf{x} \leq_{S} \mathbf{y}$ implies  $\mathbf{x} \leq_{\mathbb{N}^d} \mathbf{y}$.

\begin{remark}\label{rmk:prelim}
 {\rm    The following facts are standard and can be proved very easily.
 \begin{enumerate}
     \item If $\mathbf{x} \in \mathcal{H}(S)$ and $\mathbf{y} \leq_{S} \mathbf{x}$, then $\mathbf{y} \in \mathcal{H}(S)$.
     \item  $\mathrm{PF}(S)$ is the set of maximal elements of $\mathcal{H}(S)$ with respect to $\leq_{S}$.
     \item  For every $\mathbf{x} \in \mathcal{H}(S)$ there exists $\mathbf{f} \in \mathrm{PF}(S)$ with $\mathbf{x} \leq_{S} \mathbf{f}$.
     \item  Every maximal element of $\mathcal{H}(S)$ with respect to $\leq_{\mathbb{N}^d}$ belongs to $\mathrm{PF}(S)$.
 \end{enumerate}
 }
\end{remark}

Let $\mathrm{FA}(S)$ denote the set of all maximal elements of $\mathcal{H}(S)$ with respect to partial order $\leq_{\mathbb{N}^d}$. The set $\mathrm{FA}(S)$ is called the set of \textit{Frobenius allowable elements} of $S$~(\cite{Singhal-Lin2022}). Put $\tau(S) = \vert \mathrm{FA}(S) \vert$, then we have $\mathrm{FA}(S)\subseteq  \mathrm{PF}(S)$ and $\tau(S) \leq \operatorname{t}(S)$.

\begin{definition}
	{\rm Let $S$ be a generalized numerical semigroup and suppose $\mathcal{H}(S)$ has the maximum element, say ${\bf F}(S)$. Then, $S$ is called \textit{symmetric}  if for all ${\bf a} \in \mathbb{N}^d$: $${\bf a} \in S \iff {\bf F}(S) - {\bf a} \notin S.$$
	}
\end{definition}

 One can observe that $S$ is symmetric if and only if $\mathrm{PF}(S) = \{{\bf F}(S)\}.$ Thus, a generalized numerical semigroup $S$ is symmetric if and only if  $\tau(S) = \operatorname{t}(S) = 1.$

\section{Structure of gaps}\label{sec:gaps}

Let $d\geq 2$ and denote by $\mathbf{e}_1,\ldots,\mathbf{e}_d$ the standard basis vectors of the vector space $\mathbb{R}^d$. Let $T$ be a numerical semigroup minimally generated by the set $\{a_1,a_2, \ldots,a_n\}$. 
Given $i\in \{1,\ldots,d\}$ and  $d_j\in \mathbb{N}\setminus \{0\}$ for all $j\in \{1,\ldots,d\}\setminus \{i\}$, define the following generalized numerical semigroup in $\mathbb{N}^d$:
$$\mathcal{S}:= \mathcal{S}(T,\{d_j\}_{j\neq i})=\left\langle \{a_1 \mathbf{e}_i, a_2 \mathbf{e}_i,\ldots, a_n \mathbf{e}_i\}\cup \{\mathbf{e}_j,\mathbf{e}_i+d_j\mathbf{e}_j\mid j\in \{1,\ldots,d\}\setminus \{i\}\}\right \rangle .$$

It is easy to see that these $n+2(d-1)$ generators are minimal, so $\operatorname{e}(\mathcal{S})=e(T)+2d-2$. The importance of this family comes from the following result where it is shown that every generalized numerical semigroup of embeding dimension $2d$ is of this form. 

\begin{proposition}[\cite{Bhardwaj-Cisto2025}]\label{prop:min-emb}
Let $S\subseteq\N^d$ be a generalized numerical semigroup with $e(S)=2d$. Then, up to a permutation of coordinates, $S=\mathcal{S}(T,\{d_j\}_{j\neq i})$ for some numerical semigroup $T$ with $e(T)=2$ and some positive integers $d_j$.
\end{proposition}

Consider the numerical semigroup $T=\langle a_1, \ldots, a_n\rangle$ and $h\in \mathcal{H}(T)$, denote 
$$\alpha(h)=\max\{h-x+1\mid [x,h]\cap \mathbb{N}\subseteq \mathcal{H}(T)\}.$$
Thus $h-\alpha(h)+1,\dots,h-1,h$ are gaps of $T$, while $h-\alpha(h)\in T$. We decompose
$$\mathcal{H}(T)=I_1\sqcup I_2\sqcup\cdots\sqcup I_r$$
into maximal blocks of consecutive integers. Each block has the form $I_\nu=\{b_\nu+1,\dots,b_\nu+\ell_\nu\}$ with $b_\nu,\ b_\nu+\ell_\nu+1\in T$, and $\alpha(b_\nu+m)=m$ for $1\le m\le\ell_\nu$. 


To investigate gaps and genus of $\mathcal{S}(T,\{d_j\}_{j \neq i})$, we first investigated the genus of $\mathcal{S}(T,d_2)$ computationally in dimension two, taking $T=\langle a,b\rangle$ with $\gcd(a,b)=1$.
Using SageMath (version 10.4) \cite{sagemath}, we enumerated the gaps for several choices of $a,b$ and $d_2$, and compared the computed genus with the expression
\[
d_2\sum_{h\in \mathcal{H}(T)}\frac{\alpha(h)(\alpha(h)+1)}{2}.
\]
The computations consistently agreed with this expression and led us to examine the number of gaps with a fixed first coordinate. 
This count is established theoretically in Theorem \ref{thm:genus}, where we obtain a formula valid in arbitrary dimension and for an arbitrary numerical semigroup \(T\).
The SageMath code used for the computational experiments in this section is available at the following \href{https://github.com/priyavratd/gapsandwilf/}{GitHub repository}.  
The repository contains the notebook and instructions needed to reproduce the computations.

\begin{lemma}\label{lem:rhoproperties}
    Let $\mathbf{x} \in \mathbb{N}^d$ and define $\rho(\mathbf{x}) = \sum_{j \neq i} \lfloor \frac{x_j}{d_j} \rfloor$. Then
    \begin{enumerate}
        \item $\rho(\mathbf{x} + \mathbf{y}) \geq \rho(\mathbf{x}) + \rho(\mathbf{y})$.

        \item If $\mathbf{x} \leq_{\mathbb{N}^d} \mathbf{y}$ then $\rho(\mathbf{x}) \leq \rho(\mathbf{y})$.

        \item For $j \neq i$, $\rho(\mathbf{x} + \mathbf{e}_j)$ is either $\rho(\mathbf{x})$ or $\rho(\mathbf{x})+1$. Moreover, $\rho(\mathbf{x} + \mathbf{e}_j) = \rho(\mathbf{x})+1$ if and only if $x_j+1$ is divisible by $d_j$.
    \end{enumerate}
\end{lemma}
\begin{proof}
    The proof is easy and left to the reader.
\end{proof}

\begin{lemma}\label{lem:semigroupelements}
    Let $\mathbf{x} \in \mathbb{N}^d$. Then $\mathbf{x} \in \mathcal{S}$ if and only if $x_i - k \in T$ for some integer $k$ such that $0 \leq k \leq \rho(\mathbf{x})$.
\end{lemma}
\begin{proof}
    Note that an element $\mathbf{x}$ of $\mathcal{S}$ can be written as
    $$ \mathbf{x} = t\mathbf{e}_i + \sum_{j \neq i} \lambda_j (\mathbf{e}_i + d_j \mathbf{e}_j) + \sum_{j \neq i} \mu_j \mathbf{e}_j,$$
    where $t \in T$ and $\lambda_j, \mu_j \in \mathbb{N}$. This implies that $x_i = t + \sum_{j \neq i} \lambda_j$ and $x_j = \lambda_j d_j + \mu_j$ for $j \neq i$. Therefore, $\mathbf{x} \in \mathcal{S}$ if and only if there exist integers $\lambda_j$ such that $0 \leq \lambda_j \leq \lfloor \frac{x_j}{d_j} \rfloor$ for $j \neq i$ and $x_i -  \sum_{j \neq i} \lambda_j \in T$. Now, given such $\lambda_j$, one can take $\mu_j = x_j - \lambda_jd_j \geq 0$ and $t = x_i - \sum_{j \neq i} \lambda_j$. As $\lambda_j$ range over the intervals $[0 , \lfloor \frac{x_j}{d_j} \rfloor]$, $k = \sum_{j \neq i} \lambda_j$ takes exactly the values $0, 1, \ldots, \rho(\mathbf{x})$. This completes the proof.  
\end{proof}

\begin{theorem}\label{thm:Gaps}
     Let $\mathbf{x} \in \mathbb{N}^d$. Then $\mathbf{x} \in \mathcal{H}(\mathcal{S})$ if and only if $x_i \in \mathcal{H}(T)$ and $\rho(\mathbf{x}) \leq \alpha(x_i) - 1$. 
\end{theorem}
\begin{proof}
    Note that if $x_i \in T$ then $\mathbf{x} \in \mathcal{S}$. Suppose $x_i \in \mathcal{H}(T)$, then $x_i, x_i-1, \ldots, x_i - \alpha(x_i)+1 \in \mathcal{H}(T)$ and $x_i - \alpha(x_i) \in T$. Therefore, the least $k \geq 0$ such that $x_i - k \in T$ is $k = \alpha(x_i)$. Now the result follows by Lemma~\ref{lem:semigroupelements}. 
\end{proof}

\noindent Given a positive integer $d$ and $h \in \mathcal{H}(T)$, we also define the following subset of $\mathbb{N}^{d-1}$:
$$V_h=\left \lbrace (v_1,\ldots,v_{d-1})\in \mathbb{N}^{d-1}\mid v_i\geq 1\text{ for all }i\in \{1,\ldots,d-1\}\text{ and }\sum_{k=1}^{d-1}v_k=\alpha(h)+d-2\right \rbrace.$$
Let us denote $\{j_1< j_2 \cdots < j_{d-1}\}=\{1,\ldots, d\}\setminus \{i\}$. If $\mathbf{v}=(v_1,\ldots,v_{d-1})\in V_h$, then we define the following vector:
$$ f_\mathbf{v}=h\mathbf{e}_i+\sum_{k=1}^{d-1}(v_k d_{j_k}-1)\mathbf{e}_{j_k}.$$

The vectors $f_{\mathbf{v}}$ will be used for characterizing Frobenius allowable elements and pseudo-Frobenius elements. We can restate the Theorem~\ref{thm:Gaps} in terms of certain boxes defined with the help of $V_h$ and $f_{\mathbf{v}}$.

\begin{theorem}\label{thm:gaps}
Let $\mathbf{h} = h\mathbf{e}_i + \sum_{k=1}^{d-1}y_k \mathbf{e}_{j_k}\in \mathbb{N}^d$. Then $\mathbf{h}\in \mathcal{H}(\mathcal{S})$ if and only if  $h\in \mathcal{H}(T)$ and $\mathbf{h}\leq_{\mathbb{N}^d} f_\mathbf{v}$ for some $\mathbf{v}\in V_h$.
\end{theorem}
\begin{proof}
First suppose that $h\in \mathcal{H}(T)$ and $\mathbf{h}\leq_{\mathbb{N}^d} f_\mathbf{v}$ for some $\mathbf{v}\in V_h$. Then 
$$\rho(\mathbf{h}) \leq \rho(f_\mathbf{v}) = \sum_{k = 1}^{d-1} \left\lfloor \frac{v_k d_{j_k}-1}{d_{j_k}} \right\rfloor = \sum_{k = 1}^{d-1} (v_k-1) = \alpha(h)-1.$$
Therefore, by Theorem~\ref{thm:Gaps}, $\mathbf{h}\in \mathcal{H}(\mathcal{S})$. Conversely, suppose that $\mathbf{h}\in \mathcal{H}(\mathcal{S})$. By Theorem~\ref{thm:Gaps}, we have $h \in \mathcal{H}(T)$ and $\rho(\mathbf{h}) \leq \alpha(h) - 1$. Put $w_k = \lfloor \frac{y_{k}}{d_{j_k}} \rfloor + 1$. Therefore, $y_{k} \leq w_kd_{j_k} - 1$ and 
$$\sum_{k=1}^{d-1} w_k = \rho(\mathbf{h}) + d -1 \leq \alpha(h) + d -2.$$
Define $\mathbf{v} = (w_1+\alpha(h)-\rho(\mathbf{h})-1,w_2,\ldots, w_{d-1})$. Then $\mathbf{v}\in V_h$ and $\mathbf{h}\leq_{\mathbb{N}^d} f_\mathbf{v}$. This completes the proof. 
\end{proof}

\begin{theorem}\label{thm:genus}
	Let $\mathcal{H}(T) = I_1 \sqcup \cdots \sqcup I_r$ be the gapset of $T$ as a disjoint union of maximal consecutive gap blocks. Suppose $|I_{\nu}| = \ell_{\nu}$ for $\nu \in \{1,\ldots, r\}$. Then
	\[ \operatorname{g}(\mathcal{S}) = \operatorname{g}(\mathcal{S}(T,\{d_j\}_{j \neq i})) = \left(\prod_{j\neq i}d_j \right)
	\left( \sum_{\nu =1}^{r}
	\binom{\ell_{\nu}+d-1}{d}\right).\]
\end{theorem}
\begin{proof}
	Fix $h \in \mathcal{H}(T)$. We first prove that the number of gaps of $\mathcal{S}$ whose $i^{\mathrm{th}}$ co-ordinate is $h$ is
	\[
	\left(\prod_{j\neq i}d_j\right)
	\binom{\alpha(h)+d-2}{d-1}.
	\]
	Let $\mathbf{x}=h\mathbf{e}_i+\sum_{j\neq i}y_j\mathbf{e}_j$. Write	$y_j=q_jd_j+r_j$ with $q_j\geq 0$, and $0\leq r_j < d_j$. By Theorem~\ref{thm:Gaps}, $\mathbf{x}$ is a gap if and only if	$\sum_{j\neq i}q_j\leq \alpha(h)-1$. Therefore, we need to find all tuples $(q_1,\ldots,\widehat{q_i},\ldots,q_n)$
	such that $\sum_{j\neq i}q_j\leq \alpha(h)-1$. 	But, this number is equal to
	\[
	\binom{\alpha(h)+d-2}{d-1}.
	\]
	Also, for each fixed tuple $(q_1,\ldots,\widehat{q_i},\ldots,q_n)$, we have $d_j$ choices for the remainder $r_j$. Therefore, we have $\prod_{j\neq i}d_j$ choices for the tuples $(r_1,\ldots,\widehat{r_i},\ldots,r_n)$. Therefore, the number of gaps whose $i^{\mathrm{th}}$ co-ordinate is $h$ is
	\[ \left(\prod_{j\neq i}d_j\right)
	\binom{\alpha(h)+d-2}{d-1}.\]
	Since, the sets of gaps with different $i^{\mathrm{th}}$-co-ordinate are disjoint, we have
	\[\operatorname{g}(\mathcal{S}) = \left(\prod_{j\neq i}d_j \right) \sum_{h\in \mathcal{H}(T)} \binom{\alpha(h)+d-2}{d-1}.	\]
	Since $\mathcal{H}(T)=I_1\sqcup\cdots\sqcup I_r$. Write $I_{\nu}=\{b_{\nu}+1,\ldots,b_{\nu}+\ell_{\nu}\}$, where $b_{\nu}\in T$. In this block, one has $\alpha(b_{\nu}+m)=m$. Therefore, the block $I_{\nu}$ contributes exactly
	\[\left(	\prod_{j\neq i}d_j \right)
	\left( \sum_{m=1}^{\ell_{\nu}} \binom{m+d-2}{d-1} \right) \]
	many gaps to the genus $\operatorname{g}(\mathcal{S})$. 	Putting $k=m-1$, we have
	\[
	\left(\prod_{j\neq i}d_j\right)
	\left(
	\sum_{m=1}^{\ell_{\nu}}
	\binom{m+d-2}{d-1}
	\right)
	=
	\left(\prod_{j\neq i}d_j \right)
	\left(
	\sum_{k=0}^{\ell_{\nu}-1}
	\binom{k+d-1}{d-1}
	\right)
	= \left(\prod_{j\neq i}d_j \right)
	\binom{\ell_{\nu}+d-1}{d}. 	\] Now, by considering all the blocks $I_1,\ldots,I_r$, we get
	\[	\operatorname{g}(\mathcal{S}) = \left(\prod_{j\neq i}d_j \right)
	\sum_{\nu=1}^{r}
	\binom{\ell_{\nu}+d-1}{d}.\]
	
\end{proof}

When the gap blocks $I_{\nu}$ are known explicitly, above theorem gives very simple formulas for the genus in special cases.

\begin{corollary}
\begin{enumerate}
    \item Let $b>2$ be odd, and $T=\langle 2,b\rangle$. Then $\operatorname{g}(\mathcal{S})	= \frac{b-1}{2}
	\prod_{j\neq i}d_j.$
    \item If $d=2$ and $\{i,j\} = \{1,2\}$, then $\operatorname{g}(\mathcal{S}) =  \frac{d_j}{2}\sum_{\nu=1}^{r}\ell_{\nu}(\ell_{\nu}+1).$
    \item Let $a \geq 2$, and $T=\langle a,a+1\rangle$. Then $\operatorname{g}(\mathcal{S}) = \left( \prod_{j\neq i}d_j \right) \binom{a+d-1}{d+1}.$
\end{enumerate}
\end{corollary}
\begin{proof}
    (1) The gaps of the numerical semigroup $T$ are $\{1,3,5,\ldots,b-2\}$. Thus every block has length one, and there are $\frac{b-1}{2}$ blocks. Thus, by Theorem~\ref{thm:genus}, we get \[
	\operatorname{g}(\mathcal{S})
	= \left(\prod_{j\neq i}d_j \right)
	\left(
	\sum_{\nu =1}^{\frac{b-1}{2}}
	\binom{1+d-1}{d}
	\right)
	=
	\frac{b-1}{2}
	\prod_{j\neq i}d_j.
	\]

    (2) The result directly follows from Theorem~\ref{thm:genus}.

    (3) For $T=\langle a,a+1\rangle$, one can observe that there are $a-1$ maximal gap blocks. Moreover, for $\nu \in \{1,\ldots,a-1\}$, \[I_{\nu} = \{(\nu-1)a+\nu, \ldots, \nu a -1\}.\] 
	Therefore, by Theorem~\ref{thm:genus}, we have
	\[\operatorname{g}(\mathcal{S}) = \left( \prod_{j\neq i}d_j \right) \sum_{\nu = 1}^{a-1} \binom{\nu + d -1}{d} = \left( \prod_{j\neq i}d_j \right) \binom{a+d-1}{d+1}.\]
\end{proof}

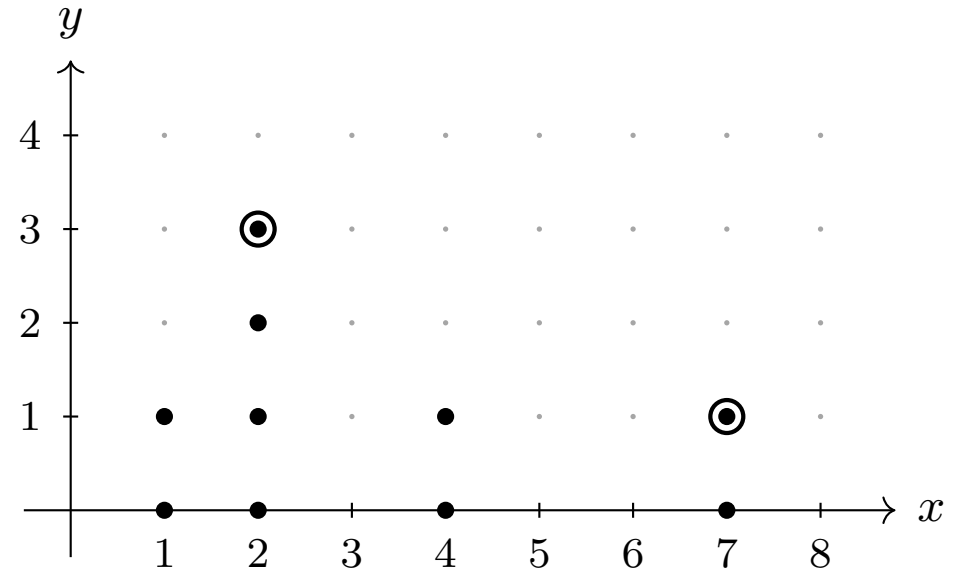
\begin{figure}[ht]
\centering
\begin{tikzpicture}[scale=0.62]

  \draw[->] (-0.5,0) -- (8.8,0) node[right] {\scriptsize $x$};
  \draw[->] (0,-0.5) -- (0,4.8) node[above] {\scriptsize $y$};
  \foreach \x in {1,...,8} \draw (\x,0.08) -- (\x,-0.08) node[below] {\tiny \x};
  \foreach \y in {1,...,4} \draw (0.08,\y) -- (-0.08,\y) node[left] {\tiny \y};
  
  \foreach \p in {(1,2), (1,3), (1,4),(2,4),(3,4),(4,4),(5,4),(6,4),(7,4),(8,4),
                  (3,3),(4,3),(5,3),(6,3),(7,3),(8,3),
                  (3,2),(4,2),(5,2),(6,2),(7,2),(8,2),
                  (8,1),(3,1), (5,1), (6,1)}
    \fill[gray!70] \p circle (0.8pt);

  \foreach \p in {(1,0),(1,1),(2,0),(2,1),(2,2),(2,3),(4,0),(4,1),(7,0),(7,1)}
    \fill \p circle (2.6pt);

  \foreach \p in {(2,3),(7,1)}
    \draw[thick] \p circle (5pt);
\end{tikzpicture}
\caption{ Example~\ref{ex:running}: gaps (filled discs), Frobenius allowable gaps (circled).}
\label{fig:running}
\end{figure}

\begin{example}\label{ex:running}
Let $d=2$, $i=1$, $T=\langle3,5\rangle$ and $d_2=2$, so that $\mathcal{S}=\langle(3,0),\allowbreak(5,0),\allowbreak(0,1),\allowbreak(1,2)\rangle$. The gaps of $T$ are $1,2,4,7$, with blocks $\{1,2\},\{4\},\{7\}$ and $\alpha(1)=1$, $\alpha(2)=2$, $\alpha(4)=\alpha(7)=1$. By Theorem~\ref{thm:Gaps}, $(h,y)$ is a gap of $S$ if and only if $h\in\{1,2,4,7\}$ and $\lfloor y/2\rfloor\le\alpha(h)-1$, so
\[
\mathcal{H}(\mathcal{S})=\{(1,0),(1,1),(2,0),(2,1),(2,2),(2,3),(4,0),(4,1),(7,0),(7,1)\}.
\]
Thus $\operatorname{g}(\mathcal{S})=10=2\big(\binom32+\binom22+\binom22\big)$, in accordance with Theorem~\ref{thm:genus}. Figure~\ref{fig:running} shows the gaps and the maximal gaps.
\end{example}

We call $h\in \mathcal{H}(T)$ a \emph{record gap} of $T$ if $\alpha(g)<\alpha(h)$ for every $g\in \mathcal{H}(T)$ with $g>h$. The Frobenius number $\mathrm{F}(T)$ is always a record gap, and a record gap is always the largest element of its block, since otherwise $h+1$ would be a larger gap with $\alpha(h+1)=\alpha(h)+1$.

\begin{theorem}\label{thm:FA}
Let $\mathbf{x}\in\mathbb{N}^d$ and $h=x_i$. Then $\mathbf{x}\in\mathrm{FA}(\mathcal{S})$ if and only if
\begin{enumerate}
\item $h$ is a record gap of $T$,
\item $x_j+1$ is divisible by $d_j$ for every $j\neq i$, 
\item $\rho(\mathbf{x})=\alpha(h)-1$.
\end{enumerate}
Equivalently, $\mathrm{FA}(\mathcal{S})=\{f_{\mathbf{v}}\mid h\text{ is a record gap of }T,\ \mathbf{v} \in V_h\}$.
\end{theorem}
\begin{proof}
Let $\mathbf{x}\in \mathcal{H}(\mathcal{S})$ such that that $h\in \mathcal{H}(T)$ and $\rho(\mathbf{x})\leq \alpha(h)-1$. We first show that $\mathbf{x}$ is not a maximal (with respect to the natural order) gap if and only if one of the following holds:
\begin{enumerate}\renewcommand{\theenumi}{\alph{enumi}}
\item $\mathbf{x}+\mathbf{e}_j\in \mathcal{H}(\mathcal{S})$ for some $j\neq i$,
\item $\rho(\mathbf{x})\leq \alpha(g)-1$ for some $g\in \mathcal{H}(T)$ with $g>h$.
\end{enumerate}
If (a) holds, $\mathbf{x}$ is clearly not a maximal gap. If (b) holds, then $\mathbf{x}+(g-h)\mathbf{e}_i$ has $i$-th coordinate $g\in \mathcal{H}(T)$. Also, $\rho(\mathbf{x}) = \rho(\mathbf{x}+(g-h)\mathbf{e}_i) \leq \alpha(g)-1$. Therefore, it is a gap by Theorem~\ref{thm:Gaps}, and thus $\mathbf{x}$ is not a maximal gap. Conversely, suppose that $\mathbf{y}>_{\mathbb{N}^d} \mathbf{x}$ is a gap. If $y_i=h$, then $y_j>x_j$ for some $j\neq i$. Therefore, $\mathbf{x}+\mathbf{e}_j\leq_{\mathbb{N}^d} \mathbf{y}$ and $\rho(\mathbf{x}+\mathbf{e}_j) \leq \rho(\mathbf{y}) \leq \alpha(h)-1$. Hence $\mathbf{x}+\mathbf{e}_j$ is a gap by Theorem~\ref{thm:Gaps} and (a) holds. If $y_i\neq h$, then $g=y_i>h$. Also, $g\in \mathcal{H}(T)$ because $y$ is a gap. Therefore, $\rho(\mathbf{x})\leq \rho(\mathbf{y})\leq \alpha(g)-1$, and (b) holds.

Now, since $\rho(\mathbf{x})\leq \alpha(g)-1$ and $\rho(\mathbf{x}+\mathbf{e}_j) \leq \rho(\mathbf{x})+1$, the vector $\mathbf{x}+\mathbf{e}_j$ is not a gap exactly when $\rho(\mathbf{x}+\mathbf{e}_j) = \rho(\mathbf{x})+1=\alpha(h)$. 
Hence (a) fails if and only if $\rho(\mathbf{x})=\alpha(h)-1$ and $x_j\equiv-1\pmod{d_j}$ for every $j\neq i$, that is, if and only if (2) and (3) hold. 
Assuming (3), condition (b) fails if and only if $\alpha(g)-1<\alpha(h)-1$ for every gap $g>h$ of $T$, that is, if and only if $h$ is a record gap. This proves the first statement. For the second statement, a vector $\mathbf{x}$ with $x_i=h$ satisfies (2) if and only if $x_{j_k}=v_kd_{j_k}-1$ for some integers $v_k \geq 1$, and then $\rho(\mathbf{x})=\sum_k(v_k-1)$. Therefore, (2) and (3) together say that $\mathbf{x}=f_{\mathbf{v}}$ with $\mathbf{v} \in V_h$.
\end{proof}

\begin{corollary}\label{cor:tau}
The number of Frobenius allowable gaps of $\mathcal{S}$ is
$$\tau(\mathcal{S})=\sum_{h\text{ record gap of }T}\binom{\alpha(h)+d-3}{d-2}.$$
\end{corollary}

\begin{proof}
Note that distinct pairs $(h,\mathbf{v})$ give distinct vectors $f_{\mathbf{v}}$. The set $V_h$ is in bijection with the compositions of $\alpha(h)+d-2$ into $d-1$ positive parts, so $|V_h|=\binom{\alpha(h)+d-3}{d-2}$. Now, the result follows from the Theorem~\ref{thm:FA}.
\end{proof}

\begin{example}\label{ex:running-FA}
Note that, in Example~\ref{ex:running}, the record gaps of $T=\langle3,5\rangle$ are $2$ and $7$. Indeed $\alpha(2)=2$ is larger than $\alpha(4)=\alpha(7)=1$, while $4$ is not a record gap since $\alpha(7)=\alpha(4)$. Hence $\mathrm{FA}(\mathcal{S})=\{(2,3),(7,1)\}$ and $\tau(\mathcal{S})=2$.
\end{example}

\section{Pseudo-Frobenius elements and type}\label{sec:PF}

The purpose of this section is to determine the set of pseudo-Frobenius elements $\mathrm{PF}(\mathcal{S})$ completely.

\begin{definition}
We call a maximal block $I$ of consecutive gaps of $T$ \emph{admissible} if $I+a_l\not\subseteq \mathcal{H}(T)$ for every $l\in\{1,\dots,n\}$.
\end{definition}

The block containing $\mathrm{F}(T)$ is always admissible. More generally, the block of a record gap is admissible. If $I=\{b+1,\dots,b+\ell\}$ and $I+a_l\subseteq \mathcal{H}(T)$, then $b+a_l\in T$ and $b+\ell+a_l$ is a gap larger than $b+\ell$ with $\alpha(b+\ell+a_l)\ge\ell=\alpha(b+\ell)$.

\begin{theorem}\label{thm:PF}
An element $\mathbf{f} \in \mathcal{H}(\mathcal{S})$ is a pseudo-Frobenius element of $\mathcal{S}$ if and only if
$$\mathbf{f}=h \mathbf{e}_i+\sum_{j\neq i}(v_jd_j-1) \mathbf{e}_j,$$
where $h$ is the largest element of an admissible block $I$ of gaps of $T$, and $(v_j)_{j\neq i}$ are positive integers with $\sum_{j\neq i}v_j=|I|+d-2$. Equivalently,
$$\mathrm{PF}(\mathcal{S})=\{f_{\mathbf{v}}\mid h\text{ is the largest element of an admissible block of }T, \mathbf{v} \in V_h\}.$$
\end{theorem}
\begin{proof}
     Let $\mathbf{f}$ be a gap. Then $h=f_i\in \mathcal{H}(T)$ and $\rho(\mathbf{f}) \leq \alpha(h)-1$. Suppose that $\mathbf{f} \in \mathrm{PF}(\mathcal{S})$. This implies that $\mathbf{f}+ \mathbf{e}_j \in \mathcal{S}$ for $j \neq i$, $\mathbf{f}+ \mathbf{e}_i +  d_j\mathbf{e}_j \in S$ for $j \neq i$, and $\mathbf{f}+ a_l \mathbf{e}_i \in \mathcal{S}$. We will explore these three cases one by one. 
 
(1) Suppose $\mathbf{f}+ \mathbf{e}_j \in \mathcal{S}$. By Theorem~\ref{thm:Gaps}, $\mathbf{f}+ \mathbf{e}_j \in \mathcal{S}$ if and only if $\rho(\mathbf{f}+ \mathbf{e}_j) \geq \alpha(h)$. Since $\rho(\mathbf{f}+ \mathbf{e}_j) \leq \rho(\mathbf{f}) +1 \leq \alpha(h)$, this happens if and only if $\rho(\mathbf{f}) =\alpha(h)-1$ and $f_j+1$ is divisible by $d_j$. Hence $\mathbf{f}+ \mathbf{e}_j \in \mathcal{S}$ for all $j\neq i$ if and only if $f_j=v_jd_j-1$ with $v_j\ge1$ for all $j\neq i$ and $\sum_{j\neq i}(v_j-1)=\alpha(h)-1$.

(2) Suppose $\mathbf{f}+ \mathbf{e}_i +  d_j\mathbf{e}_j \in \mathcal{S}$. The vector $\mathbf{f}+ \mathbf{e}_i +  d_j\mathbf{e}_j$ has $i$-th coordinate $h+1$ and satisfies $\rho(\mathbf{f}+ \mathbf{e}_i +  d_j\mathbf{e}_j) = \rho(\mathbf{f})+1 = \alpha(h)$. If $h+1\in \mathcal{H}(T)$, then $\alpha(h+1)=\alpha(h)+1>\rho(\mathbf{f}+ \mathbf{e}_i +  d_j\mathbf{e}_j)$, so it is a gap. Hence, since $d \geq 2$, these conditions hold if and only if $h+1\in T$, that is, if and only if $h$ is the largest element of its block $I$. 

(3) Assume $\mathbf{f}+ a_l \mathbf{e}_i \in \mathcal{S}$. Write $I=\{b+1,\dots,b+\ell\}$ with $b\in T$ and $h=b+\ell$, so $\rho(\mathbf{f})=\ell-1$. The vector $\mathbf{f}+ a_l \mathbf{e}_i$ has $i$-th coordinate $h+a_l$ and $\rho(\mathbf{f}+ a_l \mathbf{e}_i)=\ell-1$. If $h+a_l\in \mathcal{H}(T)$, since $b+a_l\in T$, we have $\alpha(h+a_l)\leq \ell$. Moreover, $\alpha(h+a_l)= \ell$  if and only if $b+a_l+1,\dots,b+a_l+\ell$ are all gaps.  This happens if and only if $I+a_l\subseteq \mathcal{H}(T)$. By Theorem~\ref{thm:Gaps}, $\mathbf{f}+ a_l \mathbf{e}_i \in \mathcal{S}$ if and only if $\alpha(h+a_l)\leq \ell-1$. Therefore $\mathbf{f}+ a_l \mathbf{e}_i \in \mathcal{S}$ if and only if $I+a_l\not\subseteq \mathcal{H}(T)$.

Combining the three cases, a gap $\mathbf{f}$ is pseudo-Frobenius element if and only if it has the stated form. Conversely, every element of the stated form is a gap, since $h\in \mathcal{H}(T)$ and $\rho(\mathbf{f}) = \sum_j(v_j-1)=|I|-1=\alpha(h)-1$. Finally, with $\{j_1<\cdots<j_{d-1}\}=\{1,\dots,d\}\setminus\{i\}$, the condition $\sum_jv_j=|I|+d-2=\alpha(h)+d-2$ says exactly that $\mathbf{f}=f_{\mathbf{v}}$ with $\mathbf{v}\in V_h$.
\end{proof}

\begin{corollary}\label{cor:type}
The type of $\mathcal{S}$ is
$$\operatorname{t}(\mathcal{S})=\sum_{I\text{ admissible}}\binom{|I|+d-3}{d-2}.$$
In particular, $\operatorname{t}(\mathcal{S})$ depends only on $T$ and $d$, and for $d=2$ it equals the number of admissible blocks of $T$. 
\end{corollary}

\begin{proof}
Note that distinct pairs $(h,\mathbf{v})$ give distinct vectors $f_{\mathbf{v}}$. Also $|V_h|=\binom{|I|+d-3}{d-2}$ when $h$ is the largest element of $I$. Now the result follows from the Theorem~\ref{thm:PF}.
\end{proof}

\begin{corollary}\label{cor:symmetric}
The semigroup $\mathcal{S}$ is symmetric if and only if the block of $\mathrm{F}(T)$ is the only admissible block of $T$, and either $d=2$ or $\mathrm{F}(T)-1\in T$.
\end{corollary}

\begin{proof}
As we have discussed that the block of $\mathrm{F}(T)$ is admissible. From Corollary~\ref{cor:type}, we have that $S$ is symmetric if and only if block of $\mathrm{F}(T)$ is the only admissible block of $T$ and $\binom{\alpha(\mathrm{F}(T))+d-3}{d-2} = 1$. Latter happens if and only if $d=2$ or $\alpha(\mathrm{F}(T))=1$, that is, $\mathrm{F}(T)-1\in T$. 
\end{proof}

\begin{example}\label{ex:FA-neq-PF}
Let $d=3$, $i=1$, $T=\langle4,10,17,19\rangle$, $d_2=4$ and $d_3=3$. Therefore,
$$\mathcal{S}=\langle(4,0,0),(10,0,0),(17,0,0),(19,0,0),(0,1,0),(0,0,1),(1,4,0),(1,0,3)\rangle.$$
Observe that $\mathcal{H}(T) = \{1,2,3,5,6,7,9,11,13,15\}$, and the blocks are $\{1,2,3\}$, $\{5,6,7\}$, $\{9\}$, $\{11\}$, $\{13\}$ and $\{15\}$. The blocks $\{1,2,3\}$, $\{9\}$ and $\{11\}$ are not admissible, since their translates by $4$ are $\{5,6,7\}$, $\{13\}$ and $\{15\}$. The other three blocks are admissible. For instance $10 \in \{5,6,7\}+4$, $16 \in \{5,6,7\}+10$, and $\{13\}+4=\{17\}$. By Theorem~\ref{thm:PF}, the block $\{5,6,7\}$ gives the vectors $(7,4v_2-1,3v_3-1)$ with $v_2+v_3=4$, while $\{13\}$ and $\{15\}$ give $(13,3,2)$ and $(15,3,2)$. Hence
$$
\mathrm{PF}(\mathcal{S})=\{(7,3,8),(7,7,5),(7,11,2),(13,3,2),(15,3,2)\}.$$
Now, observe that the record gaps of $T$ are $7$ and $15$. Since $\alpha(15)=\alpha(13)$, $13$ is not a record gap . Therefore, by Theorem~\ref{thm:FA}, $\mathrm{FA}(\mathcal{S})=\mathrm{PF}(\mathcal{S})\setminus\{(13,3,2)\}$. In particular, $\mathrm{FA}(\mathcal{S})\subsetneq\mathrm{PF}(\mathcal{S})$.
\end{example}

\begin{example}\label{ex:running-PF}
In Example~\ref{ex:running}, the blocks $\{1,2\}$ and $\{7\}$ of $T=\langle3,5\rangle$ are admissible, while $\{4\}$ is not, since $4+3=7$ is a gap. Hence $\mathrm{PF}(\mathcal{S})=\mathrm{FA}(\mathcal{S})=\{(2,3),(7,1)\}$ and $\operatorname{t}(\mathcal{S})=2$.
\end{example}

As a consequence of Corollary~\ref{cor:type}, the type is unbounded even in minimum embedding dimension.

\begin{corollary}\label{cor:unbounded}
Let $a\ge2$ and $T=\langle a,a+1\rangle$. Then every block of gaps of $T$ is admissible, and
$$\operatorname{t}(\mathcal{S}) =\binom{a+d-3}{d-1}.$$
In particular, for every $d\geq2$ the type of generalized numerical semigroups $S\subseteq\mathbb{N}^d$ with $e(S)=2d$ has no upper bound.
\end{corollary}

\begin{proof}
The blocks are $I_\nu=\{(\nu-1)a+\nu,\dots,\nu a-1\}$ for $1\le\nu\le a-1$. The translate $I_\nu+a$ contains $\nu a+\nu=\nu(a+1)\in T$, and $I_\nu+(a+1)$ contains $(\nu+1)a\in T$. Hence every block is admissible. Since the lengths of the blocks are $1,\dots,a-1$, Corollary~\ref{cor:type} gives
$$\operatorname{t}(\mathcal{S})=\sum_{\ell=1}^{a-1}\binom{\ell+d-3}{d-2}=\binom{a+d-3}{d-1}.$$
For $d=2$ this equals $a-1$, which is unbounded. For $d\ge3$ it is a polynomial of positive degree in $a$.
\end{proof}

  Next, we note down $d-1$ special pseudo-Frobenius elements that will play a role in Section~\ref{sec:NS}. Assume the numerical semigroup $T$ has minimal generators $\{a_1 < a_2< \cdots < a_n\}$. Let us denote $k_T=\max\{k\in \mathbb{N}\mid ka_1<a_2\}$. For all $j\in \{1,\ldots,d\}\setminus \{i\}$ denote:
 $$\mathbf{p}_j=(k_T a_1 -1)\mathbf{e}_i+((a_1-1)d_j-1)\mathbf{e}_j+\sum_{k\neq i,j}(d_k-1)\mathbf{e}_k.$$

 \begin{theorem}\label{PF:structure}
 	 The following holds:
	\begin{enumerate}
		\item $\mathbf{p}_j\in \mathrm{FA}(\mathcal{S}) \subseteq \operatorname{PF}(\mathcal{S})$ for all $j\in \{1,\ldots,d\}\setminus \{i\}$.
 		\item If $\mathbf{f}=(f_1,\ldots,f_d)\in \operatorname{PF}(\mathcal{S})$ with $\mathbf{f}\neq \mathbf{p}_j$, then $f_j<(a_1-1)d_j-1$ (that is, the $j$-th coordinate of $\mathbf{f}$ is smaller than the $j$-th coordinate of $\mathbf{p}_j$). 
 	\end{enumerate}
 \end{theorem}
 \begin{proof}
 	Observe that $h=k_Ta_1-1\in \mathcal{H}(T)$ and $\alpha(h)=a_1-1$. Therefore, one can see that $\mathbf{p}_j=f_\mathbf{v}$, with $\mathbf{v}=\alpha(h)\mathbf{e}_j+\sum_{k\neq i,j}\mathbf{e}_k$. In particular, (1) holds by Theorem~\ref{thm:FA}. In order to prove (2), it is sufficient to observe that if $\mathbf{v}\in V_h$ and $\mathbf{v}\neq \alpha(h)\mathbf{e}_j+\sum_{k\neq i,j}\mathbf{e}_k$ for all $j\in \{1,\ldots,d\}\setminus \{i\}$, then the $j$-th coordinate of $\mathbf{v}$ is strictly smaller than $\alpha(h)=a_1-1$. This concludes the proof.
 \end{proof}

\section{Generalized Wilf's conjecture}\label{sec:Wilf}
For a numerical semigroup $S,$ Wilf~(\cite{Wilf1978}) proposed a conjecture related to the Diophantine Frobenius Problem that claims that the inequality 
\[\mathrm{F}(S)+1 \leq \operatorname{e}(S) \cdot \vert \{ s \in S \mid s < \mathrm{F}(S) \} \vert \] 
is true. While this conjecture still remains open (for the state of the art of this conjecture, see~\cite{DelgadoSurvey}), two potential extensions of Wilf's conjecture to affine semigroups are studied in \cite{Cisto-et-al2020} and \cite{garcia2018extension} respectively.

\begin{conjecture}[{\bf Generalized Wilf's conjecture}{\cite{Cisto-et-al2020}}]\label{GenWilf}
Let $S\subseteq \mathbb{N}^d$ be a generalized numerical semigroup. Define:
\begin{itemize}
	\item $\operatorname{N}(S)=\{\mathbf{x}\in S\mid \mathbf{x}\leq_{\mathbb{N}^d} \mathbf{h}\text{ for some }\mathbf{h}\in  \mathcal{H}(S)\}$.
	\item $\operatorname{C}(S)=\{\mathbf{x}\in \mathbb{N}^d\mid \mathbf{x}\leq_{\mathbb{N}^d}  \mathbf{h}\text{ for some }\mathbf{h}\in  \mathcal{H}(S)\}$.
	\item $\operatorname{n}(S)=|\operatorname{N}(S)|$, $\operatorname{c}(S)=|\operatorname{C}(S)|$, and $\operatorname{g}(S) = \vert \mathcal{H}(S) \vert$.  
\end{itemize}
Then, $S$ satisfies the following inequality:
	$$\operatorname{e}(S)\operatorname{n}(S)\geq d \operatorname{c}(S)~\text{or equivalently}\quad (\operatorname{e}(S)-d)\operatorname{n}(S)\geq d \operatorname{g}(S).$$
\end{conjecture}

Let $\prec$ be a monomial order satisfying that every monomial is preceded only by a finite number of monomials. Let $Fb(S)$ denote the maximum element of $\mathcal{H}(S)$ with respect to $\prec$. By convention, $Fb(\mathbb{N}^d) = (-1,\ldots,-1)$ with $d$ coordinates. 

\begin{conjecture}[{\bf Extended Wilf's conjecture} {\cite[Conjecture 14]{garcia2018extension}}]\label{ExtWilf}
Let $S\subseteq \mathbb{N}^d$ be a generalized numerical semigroup. Let $n_{\prec}(S)$ be the cardinality of
	the finite set $\{{\bf a} \in S \mid {\bf a} \prec Fb(S)\}$, and $\mathcal{N}(Fb(S))$ be the number $n_{\prec}(S) + \operatorname{g}(S).$ Then, $S$ satisfies the following inequality:
	$$n_{\prec}(S)\operatorname{e}(S) \geq \mathcal{N}(Fb(S)) + 1.$$
\end{conjecture}

 The two conjectures are related in the context of generalized numerical semigroups in the sense that if the Conjecture~\ref{GenWilf} is true, then the Conjecture~\ref{ExtWilf} is also true (see~\cite[Proposition 6.3]{Cisto-et-al2020}). 
 
 \medskip
Observe that, since $\operatorname{e}(S)\geq 2d$ is satisfied for every generalized numerical semigroup $S\subseteq \mathbb{N}^d$, the generalized Wilf's conjecture automatically holds when $\operatorname{g}(S)\leq \operatorname{n}(S)$. Apart from equality, this fact never happens for numerical semigroups but it can be verified for generalized numerical semigroups in $\mathbb{N}^d$ with $d>1$. In particular we use this condition to prove the following general result.

\begin{theorem}\label{thm_wilf-gen}
Every generalized numerical semigroup of kind $\mathcal{S}=\mathcal{S}(T,\{d_j\}_{j\neq i})$ satisfies the generalized Wilf's conjecture.
\end{theorem}
\begin{proof}
Write $\mathcal{H}(T)=\bigsqcup_{p=1}^{r}I_p$, where 
$$I_p=\{h^{(p)}_1< h^{(p)}_2=h^{(p)}_1+1< h^{(p)}_3=h^{(p)}_1+2<\cdots < h^{(p)}_{\alpha_p}=h^{(p)}_1+\alpha_p\}\subseteq \mathcal{H}(T)$$ such that $h^{(p)}_1-1, h^{(p)}_{\alpha_p}+1\in T$ for all $p\in \{1,\ldots,r\}$. Let us denote $h^{(p)}_0=h^{(p)}_1-1\in T$ for every $p\in \{1,\ldots,r\}$. It is easy to see that $h^{(p)}_0<\operatorname{F}(T)$ and all these values are distinct considering different integers $p$. Moreover, for every $p\in \{1,\ldots, r\}$ and for every $\ell\in \{1,\ldots,\alpha_p\}$, we have $\alpha(h^{(p)}_\ell)=\ell$. By Theorem~\ref{thm:gaps}, for every $\mathbf{h}\in \mathcal{H}(\mathcal{S})$ we can express $\mathbf{h}=h_m^{(p)}\mathbf{e}_i+\sum_{k=1,k\neq i}^{d}y_k\mathbf{e}_k$ for some $p\in \{1,\ldots,r\}$ and $m\in \{1,\ldots,\alpha_p\}$, such that $0\leq y_k\leq v_kd_k-1$, where $v_k\geq 1$ for every $k\neq i$, and $\sum_{k=1,k\neq i}v_k=m+d-2$. Let us fix $j\in \{1,\ldots,d\}\setminus \{i\}$. We introduce the map $\omega:\mathcal{H}(\mathcal{S})\rightarrow \mathcal{S}$, defined by
$$\mathbf{h}=h_m^{(p)}\mathbf{e}_i+y_j\mathbf{e}_j+\sum_{\substack{k=1 \\ k\neq i,j}}^{d}y_k\mathbf{e}_k \longmapsto \omega(\mathbf{h}):=h_{\alpha_p-m}^{(p)}\mathbf{e}_i+(y_j+(\alpha_p-m)d_j)\mathbf{e}_j+\sum_{\substack{k=1 \\ k\neq i,j}}^{d}y_k\mathbf{e}_k.$$
We show that $\omega(\mathbf{h})\in \operatorname{N}(\mathcal{S})$ for every $\mathbf{h}\in \mathcal{H}(\mathcal{S})$. First observe that $\omega(\mathbf{h})\in \mathcal{S}$. This fact is trivial in the case $m=\alpha_p$, since $h^{(p)}_0\in T$. In the case $m<\alpha_{p}$, if $\omega(\mathbf{h})\in \mathcal{H}(\mathcal{S})$, by Theorem~\ref{thm:gaps} we have $y_j+(\alpha_p-m)d_j\leq w_jd_j-1$ for some integer $w_j\leq \alpha(h_{\alpha_p-m})=\alpha_p-m$, obtaining $y_j+(\alpha_p-m)d_j\leq (\alpha_p-m)d_j-1$, a contradiction since $y_j\geq 0$. So, we deduce $\omega(\mathbf{h})\in \mathcal{S}$. Moreover, let $\mathbf{w}=h_{\alpha_p}^{(p)}\mathbf{e}_i+(y_j+(\alpha_p-m)d_j)\mathbf{e}_j+\sum_{k=1,k\neq i,j}^{d}y_k\mathbf{e}_k$. We have $$y_j+(\alpha_p-m)d_j\leq (\alpha_p-m+v_j)d_j-1,$$ $y_k\leq v_kd_k-1$ for every $k\neq i,j$, and $$\alpha_p-m+v_j+\sum_{k=1,k\neq i,j}^dv_k=\alpha_p+d-2.$$ By Theorem~\ref{thm:gaps}, we deduce $\mathbf{w}\in \mathcal{H}(\mathcal{S})$. Hence, since trivially $\omega(\mathbf{h}) \leq_{\mathbb{N}^d}  \mathbf{w}$, we obtain $\omega(\mathbf{h})\in \operatorname{N}(\mathcal{S})$. Finally, we show that the map $\omega$ is injective. Let $\mathbf{h}_1,\mathbf{h}_2\in \mathcal{H}(\mathcal{S})$ such that $\omega(\mathbf{h}_1)=\omega(\mathbf{h}_2)$. In particular, we can express $\mathbf{h}_1=h_m^{(p)}\mathbf{e}_i+\sum_{k=1,k\neq i}^{d}y_k\mathbf{e}_k$ and $\mathbf{h}_2=h_n^{(q)}\mathbf{e}_i+\sum_{k=1,k\neq i}^{d}z_k\mathbf{e}_k$ for some $p,q\in \{1,\ldots,r\}$ and positive integers $y_k,z_k$. By $\omega(\mathbf{h}_1)=\omega(\mathbf{h}_2)$, we have $$h_{\alpha_p-m}^{(p)}\mathbf{e}_i+(y_j+(\alpha_p-m)d_j)\mathbf{e}_j+\sum_{k=1, k\neq i,j}^{d}y_k\mathbf{e}_k=h_{\alpha_q-n}^{(q)}\mathbf{e}_i+(z_j+(\alpha_q-n)d_j)\mathbf{e}_j+\sum_{k=1, k\neq i,j}^{d}z_k\mathbf{e}_k.$$
 From this identity we can deduce that $y_k=z_k$ for all $k\neq i,j$, $p=q$ and, as a consequence, we obtain $\alpha_p-m=\alpha_q-n$ and $y_j=z_j$. Hence $\mathbf{h}_1=\mathbf{h}_2$ and $\omega$ is injective. Therefore, we obtain $\operatorname{g}(\mathcal{S})=|\mathcal{H}(\mathcal{S})|\leq |\operatorname{N}(\mathcal{S})|=\operatorname{n}(\mathcal{S})$ and this conclude the proof.
\end{proof}


A very important consequence of Theorem~\ref{thm_wilf-gen} is that every generalized numerical semigroup having the minimum possible embedding dimension  satisfies the generalized Wilf's conjecture.

\begin{corollary}\label{cor:wilfminembdim}
Let $S\subseteq \mathbb{N}^d$ be a generalized numerical semigroup such that $\operatorname{e}(S)=2d$. Then $S$ satisfies the generalized Wilf's conjecture.
\end{corollary}
\begin{proof}
If $S\subseteq \mathbb{N}^d$ is a generalized numerical semigroup with $\operatorname{e}(S)=2d$, then $S=\mathcal{S}(T,\{d_j\}_{j\neq i})$ for some integers $d_j$ and for some numerical semigroup $T$ having $\operatorname{e}(T)=2$. Thus, the result follows by Theorem~\ref{thm_wilf-gen}.
\end{proof}


\subsection{Wilf's conjecture for generalized numerical semigroups graded by a numerical semigroup}

Let $T$ be a numerical semigroup, the $T$-graded generalized numerical semigroup $S_T^d$ in $\mathbb{N}^d$ is defined as
\[ S_T^d = \{\mathbf{x} \in \mathbb{N}^d \mid |\mathbf{x}| = x_1+\ldots+x_d \in T\}.\]
In this section we will focus on the generalized numerical semigroups graded by two generated numerical semigroups. Let $T=\langle a,b\rangle$ be a numerical semigroup. Recall that $\mathrm{F}(T)=ab-a-b$, $
\operatorname{g}(T)=\frac{(a-1)(b-1)}{2},$
and
\[
\sum_{h\in \mathcal{H}(T)}h
=
\frac{(a-1)(b-1)(2ab-a-b-1)}{12} \quad (\text{see~\cite{Brown})}.
\]

\begin{theorem}\label{thm:graded-d2}
	Let $T=\langle a,b\rangle$ be a numerical semigroup and let $S=S_T^{2}$. Then
	\[
	\operatorname{g}(S)
	=
	\frac{(a-1)(b-1)(2ab-a-b+5)}{12},
	\]
	and
	\[
	\operatorname{n}(S)
	=
	\frac{(a-1)(b-1)(4ab-5a-5b+7)}{12}.
	\]
	Moreover, $S$ satisfies the generalized Wilf's conjecture.
\end{theorem}

\begin{proof}
	By the description of the gap set of $S$ from \cite[Remark~4.2]{Cisto-Navarra}, we have
	\[
	\operatorname{g}(S)
	=
	\sum_{h\in \mathcal{H}(T)}
	\binom{h+1}{1}	=
	\sum_{h\in \mathcal{H}(T)}(h+1)
	=
	\sum_{h\in \mathcal{H}(T)}h+\operatorname{g}(T).
	\]
	Therefore,
	\begin{align*}
		\operatorname{g}(S)
		&=
		\frac{(a-1)(b-1)(2ab-a-b-1)}{12}
		+
		\frac{(a-1)(b-1)}{2}\\
		&=
		\frac{(a-1)(b-1)}{12}
		\left(2ab-a-b-1+6\right)\\
		&=
		\frac{(a-1)(b-1)(2ab-a-b+5)}{12}.
	\end{align*}
	
	Recall that $\operatorname{N}(S)
	=
	\left\{
	\mathbf{x}\in S \mid \mathbf{x} \leq_{\mathbb{N}^d}  \mathbf{h}
	\text{ for some }\mathbf{h}\in \mathcal{H}(S)
	\right\}.$ If $\mathbf{x}\in \operatorname{N}(S)$ and $s=|\mathbf{x}|$, then $s\leq \mathrm{F}(T)$. Thus, one can observe that
	\[\operatorname{N}(S) = \{\mathbf{x} \in S \mid |\mathbf{x}| \leq \mathrm{F}(T)\}.\] Moreover, since $T$ is a symmetric numerical semigroup, $s=\mathrm{F}(T)-h $ for some $h\in \mathcal{H}(T)$.  Consequently,
	\[ \operatorname{n}(S)=|\operatorname{N}(S)| = \sum_{h\in \mathcal{H}(T)}(\mathrm{F}(T)-h+1), \]
	and hence,
	\begin{align*}
		\operatorname{n}(S)
		&=
		\operatorname{g}(T)(\mathrm{F}(T)+1)
		-
		\sum_{h\in \mathcal{H}(T)}h\\
		&=
		\frac{(a-1)(b-1)(ab-a-b+1)}{2} -
		\frac{(a-1)(b-1)(2ab-a-b-1)}{12}\\
			&=
		\frac{(a-1)(b-1)}{12}
		\left(4ab-5a-5b+7\right).
	\end{align*}

	Now, it remains to prove that $S$ satisfies generalized Wilf conjecture. From \cite[Corollary~4.7]{Cisto-Navarra}, we have $\operatorname{e}(S)=a+b+2$. Since $S\subseteq\N^2$, to prove the generalized Wilf conjecture, it is sufficient to prove $(\operatorname{e}(S)-2)\operatorname{n}(S)-2 \operatorname{g}(S)\geq 0$. Put $W=(a+b)\operatorname{n}(S)-2 \operatorname{g}(S)$.
	Substituting the values of $\operatorname{g}(S)$ and $\operatorname{n}(S)$ we obtained above, we get
	\begin{align*}
		W
		&=
		\frac{(a-1)(b-1)}{12}
		\Big[
		(4ab-5a-5b+7)(a+b)-4ab+2a+2b-10 \Big]\\
		&=
		\frac{(a-1)(b-1)}{12}
		\Big[4a^2b+4ab^2-5a^2-14ab-5b^2+9a+9b-10\Big].
	\end{align*}
Put $4a^2b+4ab^2-5a^2-14ab-5b^2+9a+9b-10 = P(a,b)$. Since $b > a$, we can write $b = a+k$ with $k \geq 1$. Observe that
\[ P(a,a+k) = (4a - 5)k^2 + (12a^2 - 24a + 9)k + (8a^3 - 24a^2 + 18a - 10).\]
If $a = 2$, then $P(2,2+k) = 3k^2 +9k - 6 \geq 0$ for $k \geq 1$. If $a \geq 3$, then $4a - 5 > 0$,
$12a^2 - 24a + 9 = 12(a-1)^2 - 3 > 0$,
and
\[8a^3 - 24a^2 + 18a - 10 = 8a^2(a - 3) + 18a - 10 > 0.\]
Therefore, we get $W \geq 0$, and hence $S$ satisfies the generalized Wilf's conjecture.
\end{proof}

\begin{theorem}\label{thm:graded-d}
	Let $d\geq 3$, $3\leq a<b$, and $T$ be a numerical semigroup minimally generated by $a$ and $b$. Suppose that
	\[
	\binom{b+d-1}{d-1}
	\geq
	d\bigl(2^{d-2}(a-1)(b-1) +1\bigr).
	\]
	Then the generalized numerical semigroup $S = S_T^d	$ satisfies the generalized Wilf's conjecture. In particular, for fixed $a \geq 3$ and $d \geq 3$, the semigroup $S = S_T^d	$ satisfies the generalized Wilf's conjecture for sufficiently large $b$. 
\end{theorem}

\begin{proof}
	From \cite[Remark~4.2]{Cisto-Navarra}, we have $\operatorname{g}(S)	= \sum_{h\in \mathcal{H}(T)} \binom{h+d-1}{d-1}.$ Since $h\leq \mathrm{F}(T)$ for every $h\in \mathcal{H}(T)$, we get \[ \operatorname{g}(S) \leq \operatorname{g}(T) \binom{\mathrm{F}(T)+d-1}{d-1}.\]
	Since $\operatorname{g}(T) = 	\frac{\mathrm{F}(T)+1}{2}$, we get
	\begin{equation}\label{eq:1}
		\operatorname{g}(S) \leq \frac{\mathrm{F}(T)+1}{2}\binom{\mathrm{F}(T)+d-1}{d-1}.
\tag{1}
	\end{equation}
Let $s=\max\{ka:ka\leq \mathrm{F}(T)\}$. Then $s\in T$. If $\mathbf{x}\in \mathbb{N}^d$ satisfies $|\mathbf{x}|=s$, then $\mathbf{x}\in \operatorname{N}(S)$. Define $\mathbf{y}=\mathbf{x}+(\mathrm{F}(T)-s)\mathbf{e}_1$,	then $|\mathbf{y}| =|\mathbf{x}|+\mathrm{F}(T)-s = \mathrm{F}(T)$. Since $\mathrm{F}(T) \notin T$, we have $\mathbf{y}\in \mathcal{H}(S)$. Moreover, $\mathbf{x}\leq_{\mathbb{N}^d} \mathbf{y}$. Thus every $\mathbf{x}\in\mathbb{N}^d$ with $|\mathbf{x}|=s$ belongs to $\operatorname{N}(S)$. Consequently,
	\begin{equation}\label{eq:2}
		\operatorname{n}(S) \geq \binom{s+d-1}{d-1}. \tag{2}
	\end{equation}
Since $s$ is the largest multiple of $a$ lesser than $\mathrm{F}(T)$, we have	$0\leq \mathrm{F}(T)-s\leq a-1$. Hence $s\geq \mathrm{F}(T)-(a-1)$. Since $a\geq 3$ and $b>a$, we have $\mathrm{F}(T)=ab-a-b\geq 2(a-1)$. Therefore, $s\geq \mathrm{F}(T)-(a-1)\geq \frac{\mathrm{F}(T)}{2}$. Now, observe that
	\begin{align*}
		\frac{
			\binom{s+d-1}{d-1}
		}{
			\binom{\mathrm{F}(T)+d-1}{d-1}
		}
		&=
		\prod_{r=1}^{d-1}
		\frac{s+r}{\mathrm{F}(T)+r}.
	\end{align*}
	Since $s\geq \mathrm{F}(T)/2$, for every $r\geq1$, we have $	\frac{s+r}{\mathrm{F}(T)+r} \geq \frac{s}{\mathrm{F}(T)} \geq \frac{1}{2}$.
	Thus
	\[ \frac{\binom{s+d-1}{d-1}}{\binom{\mathrm{F}(T)+d-1}{d-1}} \geq
	\frac{1}{2^{d-1}}. \]
	Therefore, using~(\ref{eq:2}), we get
	\begin{equation}\label{eq:3}
	\operatorname{n}(S) \geq \binom{s+d-1}{d-1} \geq
		\frac{1}{2^{d-1}}
		\binom{\mathrm{F}(T)+d-1}{d-1}.
		\tag{3}
	\end{equation}
	From~\cite[Corollary~4.7]{Cisto-Navarra}, the embedding dimension of $S$ is
	\[ \operatorname{e}(S) =	\binom{a+d-1}{d-1}	+	\binom{b+d-1}{d-1},\]
	we get
	\[\operatorname{e}(S)-d \geq \binom{b+d-1}{d-1}-d.\]
	Now, using the hypothesis and~(\ref{eq:3}), we get
	\begin{align*}
		\bigl(\operatorname{e}(S)-d\bigr)\operatorname{n}(S) &\geq \left[ \binom{b+d-1}{d-1}-d \right] \frac{1}{2^{d-1}}
		\binom{\mathrm{F}(T)+d-1}{d-1}\\
		&\geq
		d\,2^{d-2}(a-1)(b-1)
		\frac{1}{2^{d-1}}
		\binom{\mathrm{F}(T)+d-1}{d-1}\\
		&=
		\frac{d(a-1)(b-1)}{2}
		\binom{F+d-1}{d-1}.
	\end{align*}
	Since $\mathrm{F}(T)+1=(a-1)(b-1)$, using~(\ref{eq:1}), we get
	\[ \frac{d(a-1)(b-1)}{2} \binom{\mathrm{F}(T)+d-1}{d-1}
	= \frac{d(\mathrm{F}(T)+1)}{2} \binom{\mathrm{F}(T)+d-1}{d-1} \geq d\operatorname{g}(S).\]
	Thus $S$ satisfies the generalized Wilf's conjecture.
\end{proof}

\section{Gaps versus non-gaps}\label{sec:gaps-nongaps}

In this section $S\subseteq\N^d$ is an arbitrary generalized numerical semigroup. For numerical semigroups, $\operatorname{n}(S)+\operatorname{g}(S)=\mathrm{F}(S)+1$ and $\operatorname{n}(S)\leq \operatorname{g}(S)$. In higher dimension, $\operatorname{n}(S)$ may exceed $\operatorname{g}(S)$ (Theorem~\ref{thm_wilf-gen}) or be much smaller (Example~\ref{ex:remove-line}). We bound $\operatorname{n}(S)$ from both sides in terms of $\operatorname{g}(S)$, the type and the number of Frobenius allowable gaps. For $\mathbf{f}\in\mathrm{PF}(S)$ let
\[
\operatorname{N}_{\mathbf{f}}=\{\mathbf{s}\in \operatorname{N}(S)\mid \mathbf{s}\leq_{\mathbb{N}^d}\mathbf{f}\},\qquad\eta_{\mathbf{f}}=|\operatorname{N}_{\mathbf{f}}|.
\]

\begin{theorem}\label{thm:g-n-upper-lower}
Let $S$ be a generalized numerical semigroup. Then
\begin{enumerate}
    \item $\operatorname{n}(S)\leq \tau(S)\operatorname{g}(S)$.
    \item $\operatorname{g}(S)\leq \sum_{\mathbf{f}\in\mathrm{PF}(S)}\eta_{\mathbf{f}} \leq \operatorname{t}(S)\operatorname{n}(S).$
\end{enumerate}

\end{theorem}

\begin{proof}
(1) Fix a monomial order $\prec$ on $\N^d$. Let $\mathbf{s}\in \operatorname{N}(S)$. Then $\mathbf{s} \leq_{\mathbb{N}^d}  \mathbf{h}$ for some gap $\mathbf{h}$, and $\mathbf{h} \leq_{\mathbb{N}^d}  \mathbf{F}$ for some $\mathbf{F}\in\mathrm{FA}(S)$, so the set $\{\mathbf{F}\in\mathrm{FA}(S)\mid \mathbf{s}\leq_{\mathbb{N}^d} \mathbf{F}\}$ is non-empty. Let $\mathbf{F}_{\mathbf{s}}$ be its least element with respect to $\prec$. The vector $\mathbf{F}_{\mathbf{s}}-\mathbf{s}$ lies in $\N^d$, and it is a gap, since otherwise $\mathbf{F}_{\mathbf{s}}=\mathbf{s} + \mathbf{F}_{\mathbf{s}}-\mathbf{s}\in S$. Define
$$\phi\colon \operatorname{N}(S)\to\mathrm{FA}(S)\times \mathcal{H}(S),\qquad\phi(\mathbf{s})=(\mathbf{F}_{\mathbf{s}},\mathbf{F}_{\mathbf{s}}-\mathbf{s}).$$
If $\phi(\mathbf{s})=\phi(\mathbf{s}')$, then $\mathbf{F}_{\mathbf{s}}=\mathbf{F}_{\mathbf{s}'}$ and $\mathbf{s}=\mathbf{F}_{\mathbf{s}}-(\mathbf{F}_{\mathbf{s}}-\mathbf{s})=\mathbf{s}'$. Hence $\phi$ is injective, and $\operatorname{n}(S)\le|\mathrm{FA}(S)|\,|\mathcal{H}(S)|=\tau(S) \operatorname{g}(S)$.

(2) Fix a monomial order $\prec$ on $\N^d$. For $\mathbf{x}\in \mathcal{H}(S)$, the set $\{\mathbf{f}\in\mathrm{PF}(S)\mid \mathbf{x}\leq_S \mathbf{f}\}$ is non-empty (see Remark~\ref{rmk:prelim}). Let $\mathbf{f}_{\mathbf{x}}$ be its least element with respect to $\prec$. Then $\mathbf{f}_{\mathbf{x}}-\mathbf{x}\in S$ and $\mathbf{f}_{\mathbf{x}}-\mathbf{x}\leq \mathbf{f}_{\mathbf{x}}$ with $\mathbf{f}_{\mathbf{x}}$ a gap, so $\mathbf{f}_{\mathbf{x}}-\mathbf{x}\in \operatorname{N}_{\mathbf{f}_{\mathbf{x}}}$. Define
\[
\psi\colon \mathcal{H}(S)\to\bigsqcup_{\mathbf{f}\in\mathrm{PF}(S)}\{\mathbf{f}\}\times \operatorname{N}_{\mathbf{f}},\qquad\psi(\mathbf{x})=(\mathbf{f}_{\mathbf{x}},\mathbf{f}_{\mathbf{x}}-\mathbf{x}).
\]
As in the proof of Part~(1), $\psi$ is injective, which gives the first inequality. The second follows from $\eta_{\mathbf{f}}\leq \operatorname{n}(S)$.
\end{proof}

The inequality $\operatorname{g}(S)\le \operatorname{t}(S)\,\operatorname{n}(S)$ was proved by Cisto and Navarra \cite{Cisto-Navarra}. The Theorem~\ref{thm:g-n-upper-lower} refines it, and is therefore a more effective tool for certifying Conjecture~\ref{GenWilf}.

\begin{example}
Let $S=\langle(2,0),(3,0),(1,1),(0,2),(0,3)\rangle$. Then
\[
\mathcal{H}(S)=\{(0,1),(1,0),(1,2),(2,1)\},\qquad \operatorname{N}(S)=\{(0,0),(0,2),(1,1),(2,0)\},
\]
and $\mathrm{PF}(S)=\{(1,2),\allowbreak(2,1)\}$, with $\operatorname{N}_{(1,2)}=\{(0,0),\allowbreak(0,2),\allowbreak(1,1)\}$ and $\operatorname{N}_{(2,1)}=\{(0,0),\allowbreak(1,1),\allowbreak(2,0)\}$. Hence $\operatorname{g}(S)=4$, while Theorem~\ref{thm:g-n-upper-lower} gives $\operatorname{g}(S)\leq 6$ and the bound $\operatorname{t}(S)\operatorname{n}(S)$ gives $\operatorname{g}(S)\le8$.
\end{example}

\begin{theorem}\label{thm:sandwich}
Every generalized numerical semigroup $S$ satisfies
\[
\frac{\operatorname{g}(S)}{\operatorname{t}(S)}\leq \operatorname{n}(S)\leq\tau(S)\,\operatorname{g}(S)\leq \operatorname{t}(S)\,\operatorname{g}(S).
\]
\end{theorem}

\begin{proof}
The first and second inequalities are from Theorem~\ref{thm:g-n-upper-lower}, and the third follows from the fact that $\tau(S) \leq \operatorname{t}(S)$.
\end{proof}

\begin{corollary}\label{cor:growth}
\begin{enumerate}
\item If $S$ is symmetric, then $\operatorname{g}(S)=\operatorname{n}(S)$.
\item For every $S$, $\operatorname{n}(S)\leq \operatorname{g}(S)^2$.
\item Let $(S_k)_{k\ge1}$ be generalized numerical semigroups with $\operatorname{t}(S_k)\le B$ for all $k$, where $B$ does not depend on $k$. Then
\[
\frac1B\leq\frac{\operatorname{n}(S_k)}{\operatorname{g}(S_k)}\le B,
\]
so $\operatorname{n}(S_k)$ and $\operatorname{g}(S_k)$ have the same order of growth. If moreover $\operatorname{g}(S_k)\to\infty$, then $\log \operatorname{n}(S_k)/\log\operatorname{g}(S_k)\to1$.
\end{enumerate}
\end{corollary}

\begin{proof}
(1) If $S$ is symmetric, $\operatorname{t}(S)=1$ and Theorem~\ref{thm:sandwich} gives $\operatorname{g}(S)\leq \operatorname{n}(S)\le \operatorname{g}(S)$. (2) Since $\mathrm{PF}(S)\subseteq \mathcal{H}(S)$, we have $\operatorname{t}(S)\leq \operatorname{g}(S)$. (3) The first claim follows from Theorem~\ref{thm:sandwich}. Taking logarithms,
\[
1-\frac{\log B}{\log \operatorname{g}(S_k)}\leq \frac{\log \operatorname{n}(S_k)}{\log \operatorname{g}(S_k)}\le1+\frac{\log B}{\log \operatorname{g}(S_k)},
\]
and both bounds tend to one when $\operatorname{g}(S_k)\to\infty$.
\end{proof}

Part (1) generalizes the classical identity $2\operatorname{g}(S)=\mathrm{F}(S)+1$ for symmetric numerical semigroups. The following example shows that the assumption on the type in part (3) cannot be removed, and that the lower bound in Theorem~\ref{thm:sandwich} is sharp.

	\begin{example}\label{ex:remove-line}
{\rm Fix $d \geq 2$, and for $k \in \mathbb{N}$, consider the sequence \[S_k = \mathbb{N}^d \setminus \{\mathbf{e}_1, 2\mathbf{e}_1, \ldots, k\mathbf{e}_1\}.\]
Observe that each $S_k$ is a generalized numerical semigroup, and \[\mathcal{H}(S_k) = \{\mathbf{e}_1, 2\mathbf{e}_1, \ldots, k\mathbf{e}_1\}.\] Therefore, $\operatorname{g}(S_k) = k$ for all $k$. Also, 
\[\operatorname{N}(S_k) = \{\mathbf{x} \in S_k \mid \mathbf{x} \leq \mathbf{h} ~\text{for some}~ \mathbf{h} \in \mathcal{H}(S_k) \}= \{\mathbf{0}\}.\] Therefore, $\operatorname{n}(S_k) = 1$ for all $k$. Moreover, for $\mathbf{f} \in \{\mathbf{e}_1, 2\mathbf{e}_1, \ldots, k\mathbf{e}_1\}$ and $\mathbf{x} \in \mathbb{N}^d \setminus \{\mathbf{0}, \mathbf{e}_1, 2\mathbf{e}_1, \ldots, k\mathbf{e}_1\}$, $\mathbf{f} + \mathbf{x} \notin \{\mathbf{e}_1, 2\mathbf{e}_1, \ldots, k\mathbf{e}_1\}$. Therefore, 
\[\mathrm{PF}(S_k) = \{\mathbf{e}_1, 2\mathbf{e}_1, \ldots, k\mathbf{e}_1\},\] and hence $\operatorname{t}(S_k) = k$ for all $k$. Therefore, the type $\operatorname{t}(S_k)$ is not bounded by a constant independent of $k$. Observe that in this example, the genus $\operatorname{g}(S_k)$ tends to infinity while $\operatorname{n}(S_k)$ remains equal to $1$. }
\end{example} 

\section{nearly-symmetric generalized numerical semigroups}\label{sec:NS}

Let $\mathbb{K}$ be a field and $R$ be a Cohen--Macaulay positively graded $\mathbb{K}$-algebra with homogeneous maximal ideal $\mathfrak{m}$. Assume that R admits a canonical module $\omega_R$ and 
$$\mathrm{tr}(\omega_R) = \sum_{\phi \in \mathrm{Hom}_{R}(\omega_R,R)} \phi(\omega_R)$$
be the trace ideal of $\omega_R$. It is known that $R$ is Gorenstein if and only if $\mathrm{tr}(\omega_R) = R.$ In \cite{Herzog-et-al2019}, Herzog et al. defined and studied nearly Gorenstein rings. A Cohen--Macaulay ring $R$ is called nearly Gorenstein if $\mathfrak{m} \subseteq \mathrm{tr}(\omega_R)$. 

\medskip
In \cite{Moscaeiello-Strazzanti2021}, Moscariello and Strazzanti studied the nearly Gorenstein property from the perspective of numerical semigroup rings. Let $S$ be a numerical semigroup minimally generated by $\{a_1, \ldots,a_n\}$. Then the semigroup ring $\mathbb{K}[S]$ is nearly Gorenstein if and only if for all $i\in\{1,\ldots,n\}$ there exists $f_i\in \operatorname{PF}(S)$ such that $a_i+f_i-f\in S$ for all $f \in \operatorname{PF}(S)$ (see \cite[Proposition~1.1]{Moscaeiello-Strazzanti2021}). This characterization is purely in terms of pseudo-Frobenius elements. As we discussed in the preliminaries that the pseudo-Frobenius elements are in one to one correspondence with the degrees of minimal generators of the last syzygy module of the semigroup ring $\mathbb{K}[S]$, then it is natural to expect that these constraints on pseudo-Frobenius elements will have some impact on the homological invariant  $\mathrm{t}(\mathbb{K}[S])$. 

\medskip
Let $\operatorname{e} \geq 4$, and $\mathfrak{S}$ be the set of all numerical semigroups of fixed embedding dimension $\operatorname{e}$. Define
$$\mathrm{t}(\mathfrak{S}) = \mathrm{max}\{\mathrm{t}(\mathbb{K}[S]) \mid S \in \mathfrak{S}\}.$$
Then it is known that $\mathrm{t}(\mathfrak{S})$ has no upper bound in terms of $\operatorname{e}$, i.e., for fixed $\operatorname{e} \geq 4$, $\mathrm{t}(\mathfrak{S})$ tends to infinity. It is an open question in the literature that when $\mathrm{t}(\mathfrak{S})$ has an upper bound in terms of $\operatorname{e}$. This question has positive answer for almost symmetric of embedding dimension lesser than $6$ (see \cite{Moscariello16, Moscariello23}), and for nearly Gorenstein numerical semigroups of embedding dimension lesser than $5$ (see \cite{Moscaeiello-Strazzanti2021}).

\medskip
The same question is also open for generalized numerical semigroups. When $S$ is an almost symmetric generalized numerical semigroup of embedding dimension $2d$ or $2d+1$, the type $\mathrm{t}(\mathbb{K}[S])$ is bounded by $3$ (see \cite{Bhardwaj-Cisto2025}). To explore boundedness of type in generalized numerical semigroups, we extend the notion of nearly Gorenstein numerical semigroups to generalized numerical semigroups.

\medskip
We extend the notion of nearly Gorenstein numerical semigroups with the help of pseudo-Frobenius elements criterion to generalized numerical semigroups, and call them nearly-symmetric generalized numerical semigroups (the reason why do we refrain to call them nearly Gorenstein is that a semigroup ring associated to a generalized numerical semigroup is never Cohen--Macaulay, hence never Gorenstein). We classify all nearly-symmetric generalized numerical semigroups of embedding dimension $2d$ and deduce that for any such nearly-symmetric generalized numerical semigroup, the type can be atmost two.

\begin{definition} \rm
	Let $S\subseteq \mathbb{N}^d$ be a generalized numerical semigroup with minimal generators $\{\mathbf{a}_1,\ldots,\mathbf{a}_n\}$. We say that $S$ is \emph{nearly-symmetric} if for all $i\in\{1,\ldots,n\}$ there exists $\mathbf{f}_i\in \operatorname{PF}(S)$ such that $\mathbf{a}_i+\mathbf{f}_i-\mathbf{f}\in S$ for all $\mathbf{f}\in \operatorname{PF}(S)$. 
\end{definition}

\begin{theorem}\label{thm:semigroup1}
	Let $S=\mathcal{S}(T,\{d_j\}_{j\neq i})\subseteq \mathbb{N}^d$ as in the described framework. If $d>2$ and $a_1\geq 3$, then $S$ is not nearly-symmetric
\end{theorem} 
\begin{proof} 
	Consider the minimal generator $\mathbf{a}=a_1\mathbf{e}_i$ and $\mathbf{f}\in \operatorname{PF}(S)$. Now, we will use the description of pseudo-Frobenius elements given in  Theorem~\ref{PF:structure}. If $\mathbf{f}\neq \mathbf{p}_j$ for every $j\in \{1,\ldots,d\}\setminus \{i\}$, then we have $\mathbf{a}+\mathbf{f}-\mathbf{p}_j\notin \mathbb{N}^d$. If $\mathbf{f}= \mathbf{p}_j$ for some $j\in \{1,\ldots,d\}\setminus \{i\}$, then we can consider $k\in \{1,\ldots,d\}\setminus \{i,j\}$, since $d>2$, and we have $\mathbf{a}+\mathbf{f}-\mathbf{p}_k\notin \mathbb{N}^d$, since $a_1\geq 3$ that implies the $k$-th coordinate is negative. This means that $S$ is not nearly-symmetric. This completes the proof.
\end{proof}

\begin{corollary}
	Let $S\subseteq \mathbb{N}^d$ be a generalized numerical semigroup with $\operatorname{e}(S)=2d$ and $d>2$. Then $S$ is nearly-symmetric if and only if $S$ is symmetric.
\end{corollary}
\begin{proof}
	If $S$ is symmetric, then it is trivially nearly-symmetric. Therefore, we may assume that $S$ is not symmetric. Observe that if $S\subseteq \mathbb{N}^d$ is a generalized numerical semigroup with $\operatorname{e}(S)=2d$, then $S=\mathcal{S}(T,\{d_j\}_{j\neq i})$ with $T=\langle a_1 , a_2 \rangle$. 
	Moreover, in this case it is known that $S$ is symmetric if and only if $a_1=2$ (see~\cite[Proposition 1.4]{Bhardwaj-Cisto2025}). Therefore, if $S$ is not symmetric then $a_1\geq 3$ and by Theorem~\ref{thm:semigroup1} we obtain $S$ is not nearly-symmetric.
\end{proof}

The case $d=2$ is shown in the following results:

\begin{lemma}\label{lemma2}
	Let $S=\mathcal{S}(T,\{d_j\}_{j\neq i})\subseteq \mathbb{N}^d$ with $d=2$ (that is, $\{i,j\}=\{1,2\}$). Then there exists a unique integer $r$ such that $\mathbf{p}=\operatorname{F}(T)\mathbf{e}_i+r\mathbf{e}_j\in \operatorname{PF}(S)$.
	In particular, if $\mathbf{f}=(f_1,f_2)\in \operatorname{PF}(S)$ with $\mathbf{f}\neq \mathbf{p}$, then $f_i<\operatorname{F}(T)$.
\end{lemma}
\begin{proof}
	Proof is immediate from Theorem~\ref{PF:structure}.
\end{proof}

\begin{theorem}
	Let $S\subseteq \mathbb{N}^2$ be a generalized numerical semigroup with $\operatorname{e}(S)=4$. Then $S$ is nearly-symmetric but not symmetric if and only if $S=\langle (3,0), (4,0), (0,1),(1,1)\rangle$ or $S=\langle (0,3), (0,4), (1,0),(1,1)\rangle$. 
\end{theorem}
\begin{proof}
	If $S=\langle (3,0), (4,0), (0,1),(1,1)\rangle$, then $\operatorname{PF}(S)=\{(2,1),(5,0)\}$ and by easy computations we obtain that $S$ is nearly-symmetric. The same for $S=\langle (0,3), (0,4), (1,0),(1,1)\rangle$ (it is obtained from the first one by permutation of coordinates). Since the two semigroups have two pseudo Frobenius elements, they are not symmetric. For the converse, suppose $S$ is nearly-symmetric not symmetric with $\operatorname{e}(S)=4$. Then $S=\langle (0,1),(1,c),(a,0),(b,0) \rangle$ or $S=\langle (1,0),(c,1),(0,a),(0,b) \rangle$ for some $a,b,c \in \mathbb{N}\setminus \{0\}$ and $\gcd(a,b)=1$. Without loss of generality we can consider the first case, and since $S$ is not symmetric we have $a\geq 3$ (see~\cite[Proposition 1.4]{Bhardwaj-Cisto2025}). Let us denote $b=ka+i$ with $k\in \mathbb{N}\setminus \{0\}$ and $i\in \{1,\ldots,a-1\}$. Let $\mathbf{p}=(ka-1,(a-1)c-1)$. 
    By Theorem~\ref{PF:structure} we have $\mathbf{p}\in \operatorname{PF}(S)$. Furthermore, let $F=\operatorname{F}(\langle a,b\rangle)=ab-a-b$ and by Lemma~\ref{lemma2} we can define $\mathbf{q}=(F,r)$ with $r\in \mathbb{N}$ such that $\mathbf{q}\in \operatorname{PF}(S)$. Since $F\neq ka-1$, we have $\mathbf{p}\neq \mathbf{q}$ and, by Theorem~\ref{PF:structure}, $r<(a-1)c-1$. Moreover, if $\mathbf{f}=(f_1,f_2)\in \operatorname{PF}(S)$ and $\mathbf{f}\neq \mathbf{q}$, then by Lemma~\ref{lemma2} we have $(0,1)+\mathbf{f}-\mathbf{q}\notin \mathbb{N}^2$. Since $S$ is nearly-symmetric, the only possibility is $(0,1)+\mathbf{q}-\mathbf{f}\in S$, which implies $r-((a-1)c-2)\in \mathbb{N}$. From this, we deduce $r=(a-1)c-2$. Furthermore, if $\mathbf{f}=(f_1,f_2)\in \operatorname{PF}(S)$ and $\mathbf{f}\notin \{\mathbf{p},\mathbf{q}\}$, then $f_1<F$ and $f_2\leq (a-1)c-2$. In particular, $(a,0)+\mathbf{f}-\mathbf{p}\notin \mathbb{N}^2$ and $(a,0)+\mathbf{q}-\mathbf{p}\notin \mathbb{N}^2$. So, the only possibility for $S$ to be nearly-symmetric is $(a,0)+\mathbf{p}-\mathbf{q}=((k+1)a-1-F,1)\in S$. As a consequence, $$(k+1)a-1-ab+a+b=(k+1)a-1-a(ka+i)+a+ka+i\geq 0.$$ Equivalently, $ka^2-(2k+2-i)a+1-i\leq 0$, which means:
	$$ \frac{2k+2-i}{2k}-\sqrt{\left( \frac{2k+2-i}{2k}\right)^2+\frac{i-1}{k}}\leq a \leq \frac{2k+2-i}{2k}+\sqrt{\left( \frac{2k+2-i}{2k}\right)^2+\frac{i-1}{k}}.$$ 
	Observe that 
	$$\frac{2k+2-i}{2k}-\sqrt{\left( \frac{2k+2-i}{2k}\right)^2+\frac{i-1}{k}}\leq 0.$$ 
	Thus the first inequality is always true. Denote $$t_{k,i}=\frac{2k+2-i}{2k}+\sqrt{\left( \frac{2k+2-i}{2k}\right)^2+\frac{i-1}{k}},$$ and observe that if $i\geq 2$, then $t_{k,i}\leq 1+\sqrt{1+\frac{1}{k}}\leq 1+\sqrt{2}<3$, obtaining $a<3$, a contradiction. Therefore, we have $i=1$, which means $a\leq \frac{2k+1}{k}$. In this case, if $k\geq 2$, then $a<3$, that is a contradiction. So, we have $k=1$ together with $i=1$, and since $a\geq 3$  we deduce $a=3$ and $b=4$. Finally, we obtained that $S=\langle (0,1),(1,c),(3,0),(4,0) \rangle$, $\mathbf{p}=(2,2c-1)$ and $\mathbf{q}=(5,2c-2)$, with $c\geq 1$. If $c\geq 2$, then $\mathbf{q}=(5,2c-2)=(4,c-2)+(1,c)\in S$, a contradiction with $\mathbf{q}\in \operatorname{PF}(S)$. Hence $c=1$ and this concludes the proof.
\end{proof}

\begin{corollary}\label{cor:NS-classification}
	Let $d$ be a positive integer and be a generalized numerical semigroup of minimum embedding dimension. Then $S$ is nearly-symmetric if and only if one of the following occurs:
	\begin{enumerate}
		\item $S$ is symmetric.
		\item $d=2$ and $S=\langle (3,0), (4,0), (0,1),(1,1)\rangle$ or $S=\langle (0,3), (0,4), (1,0),(1,1)\rangle$.
	\end{enumerate} 
\end{corollary}

We have proved in Corollary~\ref{cor:unbounded} that the type of generalized numerical semigroups $S\subseteq\N^d$ with $\operatorname{e}(S)=2d$ is unbounded. In contrast, this is not the case in nearly-symmetric generalized numerical semigroups.

\begin{corollary}\label{cor:NS-type}
	Let $S \subseteq \mathbb{N}^d$ be a nearly-symmetric generalized numerical semigroup of minimum embedding dimension then $\operatorname{t}(S) \leq 2$. 
\end{corollary}
\begin{proof}
    Follows from Corollary~\ref{cor:NS-classification}.
\end{proof}

\bibliographystyle{plain}
\bibliography{biblio-2}

	\end{document}